\input amstex
\documentstyle{amsppt}
\loadbold
\magnification = 1200
\pagewidth{5.4in}
\pageheight{7.5in}
\NoRunningHeads
\NoBlackBoxes
\define\a{\alpha}
\predefine\barunder{\b}
\redefine\b{\beta}
\define\bpage{\bigpagebreak}

\define\cd{\cdot}
\predefine\dotunder{\d}
\redefine\d{\delta}
\predefine\dotaccent{\D}
\redefine\D{\Delta}

\define\f{\frac}
\define\g{\gamma}
\define\G{\Gamma}

\define\lb\{{\left\{}
\define\la{\lambda}
\define\La{\Lambda}

\define\lra{\longrightarrow}

\define\mpage{\medpagebreak}

\define\Om{\Omega}
\define\om{\omega}
\define\op{\oplus}
\define\oper{\operatorname}

\define\ov{\overline}

\define\ox{\otimes}
\define\p{\partial}
\define\rb\}{\right\}}

\define\spage{\smallpagebreak}
\define\sbq{\subseteq}

\define\supth{{\text{th}}}
\define\T{\Theta}
\define\th{\theta}
\define\tl{\tilde}

\define\un{\underline}

\define\wh{\widehat}
\define\wt{\widetilde}
\define\x{\times}

\define\({\left(}
\define\){\right)}
\define\[{\left[}
\define\]{\right]}
\define\<{\left<}
\define\>{\right>}
\def\slantline#1#2#3#4#5{\hbox to 0pt{
}}

%SCRIPT LETTERS
%%%%%%%%%%

\def\SK{{\Cal K}}

\def\SS{{\Cal S}}
\def\ST{{\Cal T}}

\def\SW{{\Cal W}}

%%%%%%%
%%%BOLDFACE

\def\BA{{\bold A}}
\def\BB{{\bold B}}
\def\BC{{\bold C}}
\def\BD{{\bold D}}

\def\BQ{{\bold Q}}
\def\BR{{\bold R}}

\def\BZ{{\bold Z}}

\def\aut{\oper{Aut}}
\def\ext{\oper{Ext}}
\def\gl{\oper{GL}}
\def\hom{\oper{Hom}}
\def\tint{\oper{int}}
\def\img{\oper{Image}}
\def\kernel{\oper{kernel}}
\def\rank{\oper{rank}}
\def\rel{\oper{rel}}
\def\tors{\oper{Tors}}
\def\l{\big(}
\def\r{\big)}

\def\prf{\medskip\noindent{{\bf Proof.\ \ }}}
\def\pf#1{\medskip\noindent{{\bf Proof of #1.\ \ }}}
\def\rem#1{\bpage\noindent{\bf Remark #1:\ \ }}
\def\ex{\bpage\noindent{\bf Example:\ \ }}
\def\ex#1{\bpage\noindent{\bf Example #1:\ \ }}
\def\defi#1{\bpage\noindent{\bf Definition #1:\ \ }}

\topmatter
\title Dehn Surgery Equivalence Relations on Three-Manifolds
\endtitle
\author Tim D. Cochran*, Amir Gerges, and Kent Orr*\endauthor
\thanks *partially supported by the National Science Foundation
\endthanks
\affil First Version: March 28, 1997\\
This Version: January 6, 1998\endaffil
\endtopmatter

\document
\baselineskip=20pt

\subhead{\bf\S1. Introduction}\endsubhead  Suppose $M$ is an oriented
$3$-manifold. A {\it Dehn surgery\/} on $M$ (defined below) is a
process by which $M$ is altered by deleting a tubular neighborhood of
an embedded circle and replacing it again via some diffeomorphism of
the boundary torus. It was shown by W.B.R\. Lickorish [Li] and A\.
Wallace [Wa] that any closed oriented connected $3$-manifold can be
obtained from any other such manifold by a finite sequence of Dehn
surgeries. Thus under this equivalence relation all closed oriented
$3$-manifolds are equivalent. We shall investigate this same
question for more restricted classes of surgeries. In particular we
shall insist that our Dehn surgeries preserve the integral (or
rational) homology groups. Specifically, if $M_0$ and $M_1$ have
isomorphic integral (respectively rational) homology groups, is
there a sequence of Dehn surgeries, each of which preserves integral
(respectively rational) homology, that transforms $M_0$ to $M_1$?
What is the situation if we further restrict the Dehn surgeries to
preserve more of the fundamental group? Is there a difference if we
require ``integral'' surgeries? We also show that these Dehn
surgery relations are strongly connected to the following questions
concerning another point of view towards understanding $3$-manifolds.
Is there a Heegard splitting of $M_0$, $M_0=H_1\cup_fH_2$ ($H_i$ are
handlebodies of genus~$g$ and $f$ is a homeomorphism of their common
boundary surface), and a homeomorphism $g$ of $\p H_1$ such that
$M_1$ has a Heegard splitting using $g\circ f$ as the identifiation?
Since there are many natural subgroups of the mapping class group,
such as the Torelli subgroup and the ``Johnson subgroup,'' one can
ask the same question where $g$ is restricted to lie in one of these
subgroups. This is related to work of Morita on Casson's invariant
for homology $3$-spheres [Mo]. Even under these restrictions it has
been known for some time that any homology $3$-sphere is related to
$S^3$. This fact has been used to define, calculate and understand
invariants of homology $3$-spheres (such as Casson's invariant) by
choosing such a ``path to $S^3$'' in the ``space'' of $3$-manifolds.

We shall show that, in general, the equivalence relations among
$3$-manifolds thus induced are non-trivial but most can be
beautifully characterized in terms of classical invariants such as
the linking form, the cohomology ring and Massey products. This
precise and beautiful correspondence between the geometric
equivalence relation of Dehn surgery and the algebraic equivalence
of classical invariants will form the philosophical basis of a new
theory of finite type invariants for $3$-manifolds by T.D\. Cochran
and P\. Melvin [CM].

Before stating the main theorems let us be more specific. We
suppose throughout that $M$ is a compact oriented $3$-manifold.
Suppose that $\g$ is a smoothly embedded oriented circle in $M$
which is of finite order in $H_1(M;\BZ)$. Let $N(\g)$ denote a
regular neighborhood of $\g$. We define several (isotopy classes
of) embedded closed curves on $\p N(\g)$. The {\it meridian of
$\g$}, $\mu(\g)$, is given by the boundary of a meridional disk
$*\x D^2\subset S^1\x D^2\equiv N(\g)$. It is easy to see that
$\mu(\g)$ is unique and is of infinite order in $H_1(M-N(\g);\BZ)$.
A {\it parallel\/} of $\g$, denoted $\rho(\g)$, is any simple
closed curve which is homologous to $\g$ in $N(\g)$ and intersects
$\mu(\g)$ precisely once. Any two parallels ``differ'' by some
number of meridians. The longitude of $\g$, $\ell(\g)$, is a simple
closed curve which is homotopic in $N(\g)$ to some {\it positive\/}
integral multiple of $\g$ and which is of finite order in
$H_1(M-N(\g);\BZ)$. The longitude is unique and, in the case that
$\g$ is null-homologous, constitutes a {\it preferred parallel\/}.
However in general $\ell(\g)\cd\mu(\g)$ is not $1$ and hence the
longitude (unfortunately) may not be used as a parallel. With
respect to {\it some\/} choice of parallel (the preferred parallel
if $\g$ is null-homologous) one defines $p/q$-Dehn surgery on $M$
along $\g$, denoted $M_\g$, to be $[M-\dot N(\g)]\cup_\phi(S^1\x
D^2)$ where $\phi:S^1\x\p D^2\lra\p N(\g)$ is an
orientation-reversing homeomorphism which sends the curve $*\x\p
D^2$ to a curve representing $p\mu(\g)+q\rho(\g)$ in $H_1(\p N)$.
The number $p/q$ is called the {\it framing\/} of $\g$. We assume
$(p,q)=1$ and do not allow $q=0$. We call the surgery {\it
integral\/} if $q=\pm1$ and note that this is independent of the
choice of parallel. The surgery is {\it longitudinal surgery\/} if
$p\mu+q\rho$ coincides with the longitude up to sign. If $\g$ is
null-homologous this is equivalent to $p=0$. The resulting Dehn
surgery $M_\g$ does not depend on the orientation of $\g$. Then we
can ask whether or not two $3$-manifolds are related by Dehn
surgeries on curves $\g$ which are not arbitrary, but are
restricted to lie in some subgroup $N$ of $\pi_1(M)$. In addition
we shall restrict the ``framing'' $p/q$ of the surgery in order
that the integral (respectively rational) homology groups are
preserved. Specifically, suppose $N$ is the normal closure of a
finite number of elements of $G=\pi_1(M)$ and also suppose that the
elements of $N$ are trivial in $H_1(G;\BZ)$ (or $H_1(G;\BQ)$ in the
$\BQ$-case). We shall be concerned primarily with these examples:
\roster
\item"i)" $N=G_k$, the $k$-th term of the lower central series of
$G$, where $G_{k+1}=[G,G_{k-1}]$ and $G_1=G$;
\item"ii)" $N=\{e\}$, the trivial group;
\item"iii)" $N=G''$, the second derived group of $G$, that is
$G''=[G_2,G_2]$;
\item"iv)" $N=G^\BQ_2=\{x\in G\mid\exists n, x^n\in G_2\}$ which is
the set of elements which are torsion in $H_1(M;\BZ)$ (this would
be in the rational case).
\item"v)" $N=G^\BQ_k$ the $k^\supth$ term of the rational lower
central series where $G^\BQ_{k+1}$ is generated by $[G,G^\BQ_k]$
and elements of $G$ for which some power lies in $[G,G^\BQ_k]$
[Stallings].
\endroster

In all these cases $N$ is a characteristic verbal subgroup and we
can view $N$ as a {\it functor\/} from groups to normal subgroups.
That is, given a space $X$ we can speak of $N(\pi_1(X))$ in each
of these cases.

\defi{1.1} $M_1$ is $N$-surgery related to $M_0$ if there is a
finite sequence $M_0=X_0$, $X_1,\dots,X_m=M_1$ where $X_{i+1}$ is
obtained from $X_i$ by $p_i/q_i$ Dehn surgery along $\g_i$ with
$\g_i\in N(\pi_1(X_i))$ and $p_i=\pm1$ (in the $\BZ$ case where
$\g_i$ is null-homologous and $p_i$ is well-defined) and merely
non-longitudinal surgery (in the $\BQ$-case).

In particular, if $N$ is the ``$k^\supth$ lower central series
subgroup'' as in $i$), then we will say that $M_1$ is
{\it$k$-surgery equivalent to\/} $M_0$. In the next section we will
see that, in fact, the $k$-surgery relation is an equivalence
relation, and indeed preserves $\pi_1(M)/N_0$.

The case $k=2$, perhaps most fundamental, we call ``integral
homology surgery equivalence'' or sometimes merely ``surgery
equivalence''. The question of characterizing this equivalence
relation was posed by G.~Kuperberg via a newsgroup posting and was
partially answered by the Ph.D\. thesis of the second author Gerges.
A more precise answer is given in this paper. We find (see
Theorem~3.1) that this equivalence relation is completely controlled
by not just by $H_1(M)$ but by the full triple cup product
structures and by the linking form on the torsion subgroup of
$H_1(M)$. For example, we deduce that any 2 closed, connected,
oriented $3$-manifolds with isomorphic $H_1\cong\BZ^m$ ($m<3$) are
(integral homology) surgery equivalent. This generalizes the
well-known result for homology spheres ($m=0$). The latter certainly
appeared in public lectures by Andrew Casson in 1985, and we are
informed that it was known earlier. This appeared in a 1987 paper by
S.V\.~Matveev, where it was also announced that two $3$-manifolds
have isomorphic $H_1$ and linking forms if and only if they are
related by ``Borromean surgeries'' [Ma, Theorem~2 and Remark~2].

Other sample results concerning integral homology surgery
equivalence are:

\proclaim{Corollary 3.5} Let $\SS_m$ be the set of surgery
equivalence classes of closed oriented $3$-manifolds with
$H_1\cong\BZ^m$. Then there is a bijection
$\psi_*:\SS_m\to\La^3(\BZ^m)/\gl_m(\BZ)$ from $\SS_m$ to the set of
orbits of the third exterior power of $\BZ^m$ under the action
induced by $\gl_m(\BZ)$ on $\BZ^m$. Any such manifold is surgery
equivalent to one which is the result of $0$-framed surgery on a $m$
component link in $S^3$ which is obtained from the trivial link by
replacing a number of trivial $3$-string braids by a $3$-string
braid whose closure is the Borromean rings.
\endproclaim

\proclaim{Corollary 3.6} A $3$-manifold with $H_1\cong\BZ^m$ is
surgery equivalent to $\#^m_{i=1}S^1\x S^2$ if and only if its
integral triple cup product form $H^1\op H^1\op H^1\to\BZ$ vanishes
identically.
\endproclaim

\proclaim{Corollary 3.8} If $H_1(M)$ is torsion-free then $M$ is
(integral homology) surgery equivalent to $-M$.
\endproclaim

Corollary 3.8 is interesting because it is not obvious from a
geometric viewpoint how to construct such a path of surgeries, and
the result fails in general for manifolds with torsion in $H_1$.

\proclaim{Corollary 3.9} The set $\SS(\BZ_n)$ of (integral homology)
surgery equivalence classes of closed oriented $3$-manifolds with
$H_1\cong\BZ_n$ is in bijection with the set of equivalence classes
of units of $\BZ_n$, modulo squares of units, the correspondence
being given by the image of the fundamental class of the manifold
in $H_3(\BZ_n)\equiv\BZ_n$. The correspondence is also given by the
self-linking linking number of a generator of $H_1$ ($\la(1,1)=q$
and $q$ is viewed as an element of $\BZ_n$). The equivalence class
of $q\in\BZ_n$ contains the lens space $L(n,q)$.
\endproclaim

\proclaim{Corollary 3.10} ($H_1\cong\BZ_p$, $p$ prime): Any
$3$-manifold with $H_1\cong\BZ_2$ is surgery equivalent to $\BR
P(3)$. For any odd prime there are precisely two surgery
equivalence class represented by $L(p,1)$ and $L(p,q)$ where
$q\not\equiv k^2\mod p$. Hence if $p\equiv3\mod 4$ then we may
take $q=-1$ and we see that $L(p,1)$ and $-L(p,1)$ are {\it not\/}
surgery equivalent. If $p\equiv1\mod 4$ then $-L(p,1)$ is surgery
equivalent to $L(p,1)$.
\endproclaim

\proclaim{Corollary 3.12} If $\pi_1(M_0)\cong\pi_1(M_1)$ is
abelian then $M_0$ is (integral homology) surgery equivalent to
$M_1$ if and only if $M_0$ is orientation-preserving homotopy
equivalent to $M_1$.
\endproclaim

\proclaim{Proposition 3.17} Let $M_0$, $M_1$ be oriented, connected
$3$-manifolds. Then the following are equivalent:
\roster
\item"A." $M_0$ is $2$-surgery equivalent to $M_1$
\item"B." There is a framed boundary link $L$ in $M_0$ such that
Dehn surgery (with framings $\pm1$) on $L$ yields $M_1$.
\item"C." $M_0$ and $M_1$ have isomorphic linking forms and triple
cup product forms (in the sense of {\rm3.1}).
\endroster
\endproclaim

Once again the result for homology spheres, where C is vacuous,
was known. The earliest reference we can find is S.V.~Matveev
\cite{Ma; Theorem A}. This result was reproved and used by
S.~Garoufalidis \cite{Ga}.

Finally let $\SK$ be the subgroup of the mapping class group
generated by Dehn twists along simple closed curves which bound a
subsurface. Then, building on the observations of [GL] we have the
following generalization of a theorem of Morita for homology
spheres [Mo1; Proposition~2.3] [Jo1] [Mo2].

\proclaim{Theorem 3.18} Let $M_0$, $M_1$ be closed oriented
$3$-manifolds. Then the following are equivalent.
\roster
\item"A." $M_0$ and $M_1$ are $2$-surgery equivalent.
\item"B." There exist Heegard splittings $M_0=H_1\cup_fH_2$,
$M_1=H_1\cup_{\psi\circ f}H_2$ where $\psi\in\SK$.
\item"C." $M_0$ and $M_1$ have isomorphic linking forms and triple
cup product forms (as in {\rm3.1}).
\endroster
\endproclaim

The rational case for $k=2$ is controlled by $H_1(M;\BZ)/$Torsion and
the {\it integral\/} triple cup product form (Theorem~5.1). Here are
a few sample results concerning rational homology surgery
equivalence.

\proclaim{Corollary 5.2} If $m<3$ any 2 closed, oriented
$3$-manifolds with identical first Betti number $m$ are rational
homology surgery equivalent.
\endproclaim

\proclaim{Corollary 5.3} There is a bijection
$\SS^\BQ_m\to\La^3(\BZ^m)/\gl_m(\BZ)$ given by the integral triple
cup product form. Hence if $M_0$, $M_1$ have torsion-free homology
groups then they are rational homology surgery equivalent if and
only if they are integral homology surgery equivalent.
\endproclaim

The situation for higher $k$ is controlled by the above and also
higher order Massey products (see \S6).

Cases ii and iii above do not lead to equivalence relations and will
be considered in a future paper.

This paper is organized as follows:

\S1. Introduction

\S2. Preliminaries

\S3. Integral Homology Surgery Equivalence of Closed $3$-Manifolds

\S4. Proofs of Theorem 3.1 and other Basic Theorems

\S5. Rational Homology Surgery Equivalence of Closed $3$-Manifolds

\S6. Surgery Equivalence Preserving Lower Central Series Quotients

\subhead{\bf\S2. Preliminaries}\endsubhead  In this section we
prove several important technical results. In particular we show
that the surgery relations having to do with the lower central
series are in fact equivalence relations, but that symmetry fails in
general. We also show that in these cases we may safely assume our
surgeries are ``integral'' surgeries.

\proclaim{Proposition 2.1} Suppose $\g$ is a simple closed curve in
$M_0$ which is null-homologous (or merely zero in $H_1(M_0;\BQ)$ in
the $\BQ$ case). Suppose $M_1$ is the result of $\pm1/q$ surgery
along $\g$ (non-longitudinal surgery in the $\BQ$ case) and $\g'$
is the meridian of $\g$ viewed as a curve in $\pi_1(M_1)$. If
$\g\in\l\pi_1(M_0)\r_k$ for some $2\le k<\om$ then
$\g'\in\l\pi_1(M_1)\r_k$ (in the $\BQ$ case the same holds
using the rational lower central series).
\endproclaim

\pf{2.1} Let $N=N(\g)$ be the solid torus regular neighborhood of
$\g$ and $T$ its boundary. First we treat the integral case. Let
$G=\pi_1(M_0)$, $P=\pi_1(M_1)$, and suppose $\g\in G_k$. Then the
longitude $\ell(\g)$ also lies in $G_k$ and we know that there
exists an immersed $(k-1)$-stage half grope
$f:(S_1,S_2,\dots,S_{k-1})\looparrowright M_0$ whose boundary is
$\ell(\g)[FQ]$. Here we mean (as usual) that $S_1$ is a connected
surface with $f\bigm|_{\p S_1}=\g$ and that $S_i$ is a collection
of connected surfaces $S_{ij}$ such that $f\bigm|_{\p S_{ij}}$ is
$f(a_{ij})$ where $a_{ij}$ is a simple closed curve on $S_{i-1}$.
Moreover each stage $S_{i-1}$ has a $1/2$-rank system of such
curves $a_{ij}$ which occur as boundaries of $S_i$. We assume that
the immersion of $S_1$, when restricted to a collar of $S_1$, is an
embedding whose image lies in $M_0-\dot N(\g)$. It follows that
$S_1$ intersects $\g$ transversely an algebraically zero number of
times since $\ell(\g)=k\mu(\g)$ in $H_1\l M_0-\dot N(\g)\r$
has only the solution $k=0$. Recall that $\g\in G_2$ implies
$\ell(\g)$ is null-homologous in $M_0-\dot N(\g)$ and hence
$\ell(\g)\in P_2$.

By general position we also assume that all $a_{ij}$ lie in $M_0-N$
and that each $S_i$ meets $\g$ transversely. Since $\g'=\mu(\g)$ is
freely homotopic to $\pm q\ell(\g)$ in $M_1$, it follows that $\g'$
lies in $P_2$ as well. Suppose now, by induction, that
$\ell(\g)\in P_{n-1}$ for some $3\le n\le k$. We shall show that
$\ell(\g)\in P_n$ which will complete the proof in the integral
case. Note that the induction hypotheses implies that
$\mu(\g)\in P_{n-1}$. It follows that $\p S_{2j}$, for any
stage~$2$ surface, lies in $P_{n-1}$, because $\p S_{2j}$ bounds a
$(k-2)$-stage half grope in $M_0$ which lies completely in
$M_0-\tint N\sbq M_1$ except for a collection of small $2$-disks
corresponding to the transverse intersections with $\g$. Hence $\p
S_{2j}$ lies in $P_{k-1}$ modulo a product of conjugates of
$\mu(\g)$, which itself lies in $P_{n-1}$.

Now delete the algebraically-zero number of $2$-disks of
intersection of $S_1$ with $N$ and tube along $T$ to get $S^*_1$ in
$M_1-\tint N$. Then we see that $\ell(\g)$ is congruent, modulo
$P_n$, to a product of conjugates of elements of the form
$[x_i,\mu(\g)]^{\pm1}$. Since $\mu(\g)\in P_{n-1}$, $\ell(\g)\in
P_n$.

Now we address the ``rational case.'' We suppose that $\g\in
G^\BQ_k$ and hence $\ell(\g)\in G^\BQ_k$. Then there exists a
``rational'' $(k-1)$-stage half-grope whose ``boundary'' is
$\ell(\g)$. By this we mean that $\p S_1=n_1\ell(\g)$ for some
positive integer $n_1$ and similarly $\p S_{ij}=n_{ij}a_{ij}$.
Again we conclude that $S_1\cd\g$ is algebraically zero since the
equation $n_1\ell(\g)=m\mu(\g)$ in $H_1(M_0-\dot N(\g);\BZ)$ has
only the solution $m=0$ since $\ell(\g)$ is torsion while $\mu(\g)$
is not. We claim that $\g'=\mu(\g)$ is ``rationally related'' to
$\ell(\g)$ in $P$, that is that there exist integers $x$, $y$ such
that $\l\mu(\g)\r^x=\l\ell(\g)\r^y$ where $x\neq0$. To
see this, consider the (abelian) subgroup $T$ of $P$ generated by
$\mu$ and the parallel $\rho$. Suppose $\ell(\g)=a\mu+b\rho$ where
$(a,b)=1$. Then $T$ is a quotient of the abelian group
$A=\<\mu,\rho\mid p\mu+q\rho=0\>$ and contains $\ell(\g)$. The
vectors $(p,q)$ and $(a,b)$ are linearly independent in $\BZ\x\BZ$
since they are primitive and $(p,q)\neq\pm(a,b)$ since our surgery
is not longitudinal. Hence $A/\<\ell(\g)\>$ is a finite group and
thus there are integers $x$, $y$ $x\neq0$ such that
$x\mu=y\ell(\g)$ in $A$ and hence in $T$. Using this relation, the
proof now proceeds as in the integral case. \qed

\proclaim{Corollary 2.2} The relation of $k$-surgery equivalence is
an equivalence relation on the set of oriented $3$-manifolds. The
relation of rational $k$-surgery equivalence is also an equivalence
relation (Here we mean the rational lower central series and
non-longitudinal surgery as in {\rm example v)}.
\endproclaim

\pf{2.2} Reflexivity and transitivity are obvious and symmetry is
guaranteed by 2.1.

\proclaim{Proposition 2.3} If $M_0$ and $M_1$ are $k$-surgery
equivalent (respectively rationally $k$-surgery equivalent) then
they are so equivalent using only {\bf integral} surgeries, that is
$\pm1$ surgeries (respectively integral non-longitudinal
surgeries). In generality, if $M_1$ is $N$-surgery related to $M_0$
then there is an epimorphism
$\pi_1(M_1)\twoheadrightarrow\pi_1(M_0)/N\l\pi_1(M_0)\r$ (in
either $\BZ$ or $\BQ$ case of {\rm Definition~1.1}). Consequently if
$M_1$ is rationally $2$-surgery equivalent to $M_0$ then $\b_1(M_0)=
\b_1(M_1)$. If $M_1$ is integrally $2$-surgery equivalent to $M_0$
then $H_1(M_0;\BZ)\cong H_1(M_1;\BZ)$.
\endproclaim

\pf{2.3} The sequence of homeomorphisms shown in Figure~2.4 using
the ``Rolfsen-Kirby'' calculus is well-known (see [CG; p\.~501] [R;
p\.~ ]). This shows that $1/n$ surgery on $\g$ is the same as a
sequence of $\pm1$ surgeries on parallel copies of $\g$, denoted
$\g_1$, $\g_2,\dots,\g_n$. Let $M^*_i$ be the result of surgery on
$\{\g_1,\g_2,\dots,\g_i\}$. We may view $\g_{i+1}$ as the longitude
of $\g_i$ and 2.1 guarantees that if $\g_i$ lies in the $k$-th term
of the lower central series of $\pi_1(M^*_{i-1})$ then $\mu(\g_i)$
and $\ell(\g_i)$ lie in $\l(\pi_1(M^*_i)\r_k$. Hence $M^*_n$
is $k$-surgery equivalent to $M_0$ via $+1$ surgeries as claimed.

\midinsert
\vspace{1.2in}
\botcaption{Figure 2.4}\endcaption
\endinsert

Now consider the case that $M_1$ is rationally $k$-surgery
equivalent to $M_0$ via a single surgery on $\g$ in $M_0$. Since
$M_0=S^3_J$ for some framed link $J$ in $S^3$, which may be assumed
to be disjoint from $\g$, we may consider $\g\subset S^3$ with
framing $p/q$ with respect to the longitude of $\g$ in $S^3$.
Suppose the surgery is {\it not\/} integral, i.e., $q\neq1$ (we
assume $q>0$). Then $p/q=\pm\l m+\f rq\r=
\pm\l(m+1)-\f{q-r}q\r$ where $m=\[|p/q|\]$ and $0<r<q$.
Then it is well known that the 3 pictures of Figure~2.5 are
homeomorphic [Rolfsen,\ \ \ \ ], where the upper sign is used if
$p\ge0$. Here the framings are all relative to  $S^3$. Since
$0<r<q$ and $0<q-r<q$, we may use 2.5 as the basis of an induction
to reduce all surgeries in a daisy-chain of circles to integers
(get $q=1$). Since the first circle $\g$ lies in
$\l\pi_1(M_0)\r^\BQ_k$, it suffices to show that the second
curve in 2.5, say $\g_2$, lies in $G^\BQ_k$ where $G=\pi_1(M^*)$,
$M^*$ being the result of $\pm m$ or $\pm(m+1)$ surgery on $\g$. In
addition we must show that the $\pm m$ (or $\pm(m+1)$) surgery is
integral and non-longitudinal and show the surgery on $\g_2$ is
also non-longitudinal. Firstly, {\it at most one\/} of $\pm m$ or
$\pm(m+1)$ could be longitudinal with respect to $M_0$ so we may
choose the non-longitudinal one. Moreover the surgery on $\g$ is
integral means that the defining torus homeomorphism extends over a
solid torus. But this is independent of coordinate system so the
fact that $\pm m$, $\pm(m+1)$ are integers suffices to show these
surgeries are integral relative to $M^*$.

\midinsert
\vspace{1.2in}
\botcaption{Figure 2.5}\endcaption
\endinsert

Now, as has been mentioned, and is an immediate consequence of the
second part of 2.3 (which is proved below), rational $k$-surgery
equivalence preserves $H_1/$(torsion). Thus
$\b_1(M_0)=\b_1(M^*)=\b_1(M_1)$. This implies that the framing on
$\g_2$ relative to $M^*$ is non-longitudinal since it is precisely
the longitudinal surgery which changes $\b_1$.

In general, if $M_1$ is $N$-surgery related to $M_0$ via a single
$\pm1/n$ surgery then Figure~2.4 shows that there is a link $L$ of
$n$ components in $M_0$, each component of which lies in
$N\l\pi_1(M_0)\r$, such that adding $n$ two-handles to
$M_0\x[0,1]$ along $L$ yields a cobordism $W$ between $M_0$ and
$M_1$ rel~$\p M_0=\p M_1$. But then both inclusion maps induce
epimorphisms on $\pi_1$ and the kernel of $\pi_1(M_0)\to\pi_1(W)$
lies in $N$. Thus
$\pi_1(W)/N\l\pi_1(W)\r\cong\pi_1(M_0)/N\l\pi_1(M_0)\r$
and the second claimed result follows easily in the $\BZ$ case. In
the $\BQ$-case, if $M_1$ is $N$-related to $M_0$ via a single
non-longitudinal surgery then 2.5 shows that there is a cobordism
$W$ as above and a link $L$ each of whose components lies in
$N\l\pi_1(M_0)\r$ (in fact most are null-homotopic!). Then
the argument above for the $\BZ$-case holds. Note that if $M_1$ is
rationally $2$-surgery equivalent to $M_0$ then $\pi_1(M_1)$ maps
onto the free abelian group
$\BZ^{\b_1(M_0)}=\pi_1(M_0)/(\pi_1(M_0))^\BQ_2$. Hence
$\b_1(M_1)\ge\b_1(M_0)$. But by symmetry (2.1)
$\b_1(M_0)=\b_1(M_1)$ and consequently
$\pi_1(M_1)/(\pi_1(M_1))^\BQ_2\cong\BZ^{\b_1(M_0)}$. \qed

\midinsert
\vspace{1.2in}
\botcaption{Figure 2.7}\endcaption
\endinsert

\ex{2.6} It is important to note that ``$N$-surgery related'' is
{\it not\/} a symmetric relation if $N=\{e\}$ or $N=G''$. In a
later paper we will consider strengthening the $N$-surgery
relation to {\it force\/} symmetry. Figure~2.7a shows a ``Kirby
calculus'' description of $M_0=S^1\x S^2$ with a dashed curve $\g$
which is clearly null-homotopic in $M_0$. Yet $+1$ surgery along
$\g$ yields $M_1$ which (since the Whitehead link is symmetric) is
homeomorphic to the manifolds of Figure~2.7b and 2.7c. Hence $M_1$
is $0$-surgery on a left-handed trefoil knot and the loop $\g'$ is
neither null-homotopic in $M_1$ nor even in
$\l\pi_1(M_1)\r''$. If $M_0$ were $N$-surgery related to
$M_1$ for $N=G''$ (or $\{e\}$) then by 2.3 there would exist an
epimorphism from $\BZ=\pi_1(M_0)$ to
$\pi_1(M_1)/N\l\pi_1(M_1)\r$. Since the trefoil knot has
non-trivial Alexander module, this is not possible. One also sees
that, in the case $N=\{e\}$, forcing symmetry would force $\pi_1$
itself to be preserved (since $3$-manifold groups are Hopfian) and
{\it perhaps\/} this is too strong to be of interest.

We have defined our ``equivalence'' relations to be generated by
single Dehn surgeries. It is also possible to define relations
generated by surgeries on certain types of {\it links\/}. These are
sometimes equivalent notions as the following show. The proof of the
first is elementary and left to the reader.

\proclaim{Proposition 2.8} The following are equivalent.
\roster
\item"1." $M_0$ and $M_1$ are $2$-surgery equivalent.
\item"2." There is a link $L=\{L_1,\dots,L_m\}$ in $M_0$ with
null-homologous components, each framed $\pm1$, with $\ell
k(L_i,L_j)=0$ so that $M_1$ is the result of surgery along $L$. In
other words, the ``linking matrix'' of $L$ is invertible over $\BZ$
and diagonal.
\endroster
\endproclaim

\proclaim{Proposition 2.9} The following are equivalent.
\roster
\item"1." $M_0$ and $M_1$ are rationally $2$-surgery equivalent.
\item"2." There is a framed link $L=\{L_1,\dots,L_m\}$ in $M_0$, each
component of which is rationally null-homologous in $M_0$, such that
the ``linking matrix'' of $L$ is non-singular over $\BQ$ and such
that $M_1$ is obtained by surgery on $L$.
\endroster
\endproclaim

Here, by ``linking matrix'' of the framed link
$L=\{\g_1,\dots,\g_m\}$ we mean the matrix over $\BQ$ given by
$v_{ij}=\ell k(\rho_i,\g_j)$ where $\rho_i$ here is the circle on
$\p N(\g_i)$ which bounds the meridional disk in the surgery solid
torus ($p_i\mu_i+q_i\rho_i$ in the notation of \S1). The proof of
2.9 is given after the proof of Theorem~4.2.

\subhead{\bf\S3. Integral Homology Surgery Equivalence of Closed
3-manifolds}\endsubhead  In this chapter we give a comprehensive
treatment of $2$-surgery equivalence of closed $3$-manifolds. This
is precisely the equivalence relation generated by Dehn surgeries
which preserve integral homology and thus was called {\it
HTS-equivalence\/} (homologically trivial surgery) by Kuperberg
[Ku] and Gerges [Ge]. This equivalence relation is perhaps the most
basic and important. It forms the basis of the philosophy of
Cochran and P\.~Melvin in their theory of finite type invariants
for arbitrary $3$-manifolds [CM]. The question of characterizing
$2$-equivalence was asked by Kuperberg and answered by Gerges in
his Ph.D\. thesis. Here we prove a sharper theorem. Our
characterization theorem says that $M_0$ and $M_1$ are HTS
equivalent precisely when they have the same $H_1$ and some
isomorphism induces an isomorphism of $\BQ/\BZ$ linking forms and
that part of the cohomology ring coming from triple cup products. In
the next chapter we will prove the characterization theorem, which
appears in Gerges [Ge] without the relation to the linking form. In
this chapter we will discuss examples, invariants and
representatives for the $2$-equivalence classes.

Before stating the theorem, we set up some notation. We let
$K(H_1(M_0),1)$ be the usual Eilenberg-Maclane space with
fundamental group $H_1(M_0)$. We can build this space from $M_0$ by
adding cells of dimension greater than $1$ and we let
$f_0:M_0\to K(H_1(M_0),1)$ denote this inclusion. Then if
$\phi_1:H_1(M_1)\to H_1(M_0)$ is any isomorphism, there is a unique
homotopy class $f_1:M_1\to K(H_1(M_0),1)$ inducing $\phi_1$ on
$H_1$. Let $B: H^1(\un{\ \ };\BZ_n)\to H^2(\un{\ \ };\BZ)$ denote the
Bockstein operator associated with the short exact sequence
$0\lra\BZ\overset n\to\lra\BZ\overset\tau\to\lra\BZ_n\lra0$.

The following theorem in the case of homology $3$-spheres was
certainly known to and used by Andrew Casson in public lectures at
M.S.R.I\. in 1985. We are informed that it was known even earlier. In
this case it says merely that any two oriented homology $3$-spheres
are $2$-surgery equivalent. This case also appeared in a 1987 paper
of S.V\.~Matveev. In the latter, moreover, it is proved that two
$3$-manifolds have isomorphic $H_1$ and linking forms if and only if
they are related by ``Borromean surgeries,'' a result clearly close
in spirit to our final one [Ma, Theorem~2 and Remark~2].

The equivalence of B and (a slightly stronger version of) D is
claimed in passing in [Tu1], but no proof is offered.

\proclaim{Theorem 3.1} {\rm (see [Ge])} Suppose $M_0$ and $M_1$ are
closed, oriented, connected $3$-manifolds. The following 4
conditions are equivalent.

\mpage
\noindent{\rm A)} $M_0$ and $M_1$ are $2$-surgery equivalent,
i.e., each can be obtained from the other by a sequence of $\pm1$
surgeries (equivalently $\pm1/q$ surgeries) or null-homologous
circles.

\spage
\noindent{\rm B)} There exists an isomorphism $\phi_1:H_1(M_1)\to
H_1(M_0)$ such that $(f_0)_*([M_0])=(f_1)_*([M_1])$ in
$H_3(H_1(M_0);\BZ)$ where $f_0$, $f_1$ are as above, In brief one
could also say that $M_0$ and $M_1$ have the same homology and are
{\it bordant\/} over $K(H_1(M_0),1)$ for some $f_i$.

\spage
\noindent{\rm C)} There exists an isomorphism $\phi_1:H_1(M_1)\to
H_1(M_0)$ such that the set of induced maps
$\phi^1_n:H^1(M_0;\BZ/n\BZ)\to H^1(M_1;\BZ/n\BZ)$ for $n=0$ and for
each $n=p^r$ where $p^r$ is the exponent of the $p$-torsion subgroup
of $H_1(M_0;\BZ)$ (all elements of order some power of the prime
$p$) satisfies the following:
\roster
\item"a)"
$\<\a\cup\b\cup\g,[M_0]\>=\<\phi^1_n(\a)\cup\phi^1_n(\b)\cup
\phi^1(\g),[M_1]\>$ where $\a$, $\b$, $\g\in H^1(M_0;\BZ/n\BZ)$ and
$[M_i]$ denotes the fundamental class in $H_3(M_i;\BZ/n\BZ)$,

\item"b)" $\<\a\cup\tau_* B(\g),[M_0]\>=\<\phi^1_n(\a)\cup\tau_*
B\(\phi^1_n(\g)\),[M_1]\>$ where $\a$, $\g$, $B$ are as above, but
$n\neq0$, and $\tau:H^2(M_i;\BZ)\to H^2(M_i;\BZ_n)$.
\endroster

\noindent{\rm  D)} The same condition as {\rm C} with {\rm b)}
replaced by
\roster
\item"c)" If $\la_i$ represent the $\BQ/\BZ$ linking forms on
$T(H_1(M_i))$ then $\la_1(x,y)=\la_0(\phi_1(x),\phi_1(y))$ for all
$x$, $y\in T(H_1(M_1))$, that is to say that $\phi_1$ induces an
isomorphism between $\la_0$ and $\la_1$.
\endroster
\endproclaim

Let $A$ be a finitely generated abelian group and
$A^*_n=\hom(A;\BZ_n)\equiv H^1(A;\BZ_n)$ for $n=0$ or $n=p^r$ (the
exponent of the $p$-torsion subgroup of $A$). Consider a set of
skew-symmetric trilinear forms $u_n:A^*_n\x A^*_n\x
A^*_n\lra\BZ_n$, where $n$ ranges over $\{0,p^r\}$ as above, which
are compatible in the sense of [Tu2; Definition~1.2]. Let $\la:$
Torsion~$A\x$ Torsion~$A\lra\BQ/BZ$ be a non-degenerate symmetric
bilinear form. Any automorphism $\phi:A\to A$ induces isomorphic
forms $\{\phi^*(u_n)\}$ and $\{\phi_*\la\}$ given by
$\phi^*(u_n)(\a,\b,\g)=u_n\(\phi^*_n(\a),\phi^*_n(\b),\phi^*_n(\g)\)$
where $\phi^*_n:A^*_n\to A^*_n$ and
$\phi_*\la(x,y)=\la(\phi^{-1}x,\phi^{-1}y)$. Given any oriented
$3$-manifold with $H_1\cong A$, the triple cup product forms and
linking form yield a pair $(\{u_n\},\la)$ which is well-defined up
to isomorphism. Let $\SS(A)$ be the set of isomorphism classes of
such pairs which are realizable by a $3$-manifold. In fact by [Tu2;
Theorem~1] and [KK], {\it any\/} pair is realizable if $A$ has no
$2$-torsion. In general there is a mild compatibility condition
between $u_{2^r}$ and $\la$. Hence Theorem~3.1 may be restated as
follows.

\proclaim{Theorem 3.1 (Restatement)} Let $A$ be a finitely generated
abelian group. The set of surgery equivalence classes of closed,
oriented $3$-manifolds with $H_1\cong A$ is in bijection with
$\SS(A)$.
\endproclaim

\ex{3.2 $(H_1\cong\BZ^m\ m<3)$} Any two homology $3$-spheres are
surgery equivalent since $B$ is trivially satisfied in this case.
Indeed since $H_3(\BZ^m)=0$ if $m<3$, any two $3$-manifolds $M_0$,
$M_1$ with $H_1\cong\BZ^m$ are surgery equivalent if $m<3$.

\ex{3.3 $H_1\cong\BZ^3$} If $M_0=\#^3_{i=1}S^1\x S^2$ and
$M_1=S^1\x S^1\x S^1$ then $M_0$ is not surgery equivalent to $M_1$
because the image of $M_0$ in $H_3\l H_1(M_0)\r$ is zero since it
factors through $H_3\l\pi_1(M_0)\r$. But for {\it any\/}
automorphism of $\BZ\x\BZ\x\BZ$, the induced map $M_1\to S^1\x S^1\x
S^1=K\l H_1(M_0),1\r$ is of degree $\pm1$ since the identity map is
degree~1 and Aut$\l H_3(\BZ^3)\r=\{\pm1\}$. Hence condition~B fails
(equivalently condition~C part~a with $n=0$). More generally, let
$M_n$ be the $3$-manifold shown in figure~3.4 as zero surgery on a
link with $\ov\mu(123)=n$. Then $M_n$ is surgery equivalent to
$M_m$ if and only if $|n|=|m|$ since the triple cup product of the
Hom-duals of the meridians is $n$ times the fundamental class in
$H^3$. Any element of Aut$(\BZ^3)$ induces an element of Aut$\l
H_3(\BZ^3)\r$ and this correspondence is $P\to\det P$. Since $\det
P=\pm1$, the classes $n$, $m$ in $H_3(\BZ^3)$ are equivalent under
the action of $\gl(3,\BZ)$ if and only if $n=\pm m$. Note that this
implies that $M_0$ and $M_1$ are surgery equivalent if and only if
they have identical Lescop invariant [Les]. This generalizes to give
the following.

\midinsert
\vspace{1.2in}
\botcaption{Figure 3.4}\endcaption
\endinsert

\proclaim{Corollary 3.5} Let $\SS_m$ be the set of surgery
equivalence classes of closed oriented $3$-manifolds with
$H_1\cong\BZ^m$. Then there is a bijection
$\psi_*:\SS_m\to\La^3(\BZ^m)/\gl_m(\BZ)$ from $\SS_m$ to the set of
orbits of the third exterior power of $\BZ^m$ under the action
induced by $\gl_m(\BZ)$ on $\BZ^m$. Any such manifold is surgery
equivalent to one which is the result of $0$-framed surgery on a $m$
component link in $S^3$ which is obtained from the trivial link by
replacing a number of trivial $3$-string braids by a $3$-string
braid whose closure is the Borromean rings.
\endproclaim

\pf{3.5} Fix an identification $H_3\l(S^1)^m\r\equiv
H_3(\BZ^m)\equiv\La^3\BZ^m$. For any $M$ with $H_1\cong\BZ^m$
choose an isomorphism $\phi:H_1(M)\to\BZ^m$. This induces a unique
homotopy class of maps $\psi:M\to(S^1)^m$. The image of the
fundamental class of $M$ is the desired element $\psi_*([M])$. All
possible isomorphisms $\phi$ may be achieved by post-composing a
fixed $\phi$ with elements of $\gl_m(\BZ)$. Then $\psi_*(M)$ is
well-defined in the orbit space. If $M$ and $M'$ are surgery
equivalent then condition~B of 3.1 guarantees their images are the
same. Hence $\psi_*$ is well-defined. $B\Rightarrow A$ implies that
$\psi_*$ is injective. The surjectivity of $\psi_*$ follows from
work of D\.~Sullivan [Su]. Alternatively, for any set of
$\pmatrix m\\ 3\endpmatrix$ integers $\{a_{ijk}\mid1\le i<j<k\le m\}$
is it easy to construct an ordered link of $m$ components in $S^3$
such that $\ov\mu(ijk)=a_{ijk}$ by the procedure described in the
last sentence of the Corollary. If $M$ is zero surgery on this link
then $\psi_*(M)=\sum a_{ijk}(e_i\land e_j\land e_k)$ with respect to
a basis induced by the meridians (see Lemma~4.2 of
\cite{Tu2}). \qed

\bpage
The structure of the set $\La^3(\BZ^m)/\gl_m$ seems to be quite
complicated for large $m$ and so the general decidability question
for whether or not two $3$-manifolds are surgery equivalent may not
be easy. However since the $0$ element is the only element in its
orbit we have:

\proclaim{Corollary 3.6} A $3$-manifold with $H_1\cong\BZ^m$ is
surgery equivalent to $\#^m_{i=1}S^1\x S^2$ if and only if its
integral triple cup product form $H^1\op H^1\op H^1\to\BZ$ vanishes
identically.
\endproclaim

We also observe the following surprising result.

\proclaim{Corollary 3.7} The map $\SS_3\overset f\to\lra\SS_4$
given by $M\to M\#S^1\x S^2$ is a bijection. Thus any $N$ with
$H_1(N)\cong\BZ^4$ is surgery equivalent to precisely one of
$M_n\#S^1\x S^2$ where $n\ge0$ (see Figure~3.4).
\endproclaim

\pf{3.7} $\La^3\BZ^4\cong\La^1\BZ^4\cong\BZ^4$ by duality. Thus
$\La^3\BZ^4/\gl_4\cong\BZ^4/\gl_4\cong\BZ_+\cup\{0\}$. Under this
bijection, $n(e_1\land e_2\land e_3)$ goes to $ne_4$ and the former
is $\psi_*(M_n\#S^1\x S^2)$. \qed

\proclaim{Corollary 3.8} If $H_1(M)$ is torsion-free then $M$ is
surgery equivalent to $-M$.
\endproclaim

\pf{3.8} The element of $\gl_m$ which reverses the order of a basis
$\{e_1,\dots,e_m\}\to\{e_m,\dots,e_1\}$ induces $-1$ on
$\La^3\BZ^m$. \qed

\bpage
Now consider that $H_1\cong\BZ_n$. First consider the general
question of which classes $\mu\in H_3(A)$ can be realized as the
image of the fundamental class of a $3$-manifold $M_0$ under some
map $\pi_1(M_0)\twoheadrightarrow H_1(M_0)\overset f\to\lra A$
where $f$ is an {\it isomorphism\/}. This question has been
answered by Turaev in great generality. The answer is that $\mu$ is
realizable if and only if $x\mapsto x\cap\mu$ is an isomorphism
$\tors H^2(A)\to\tors H_1(A)$. Both groups are $\BZ_n$ in the case
at hand. If we denote by $\mu=1\in H_3(\BZ_n)=\BZ_n$ the image of
the class of $L(n,1)$ under some map then certainly $x\mapsto
x\cap1$ is an isomorphism. A general class $\mu=k\cd1$ will induce
the map $x\mapsto k(x\cap1)$ which is the composition of the
isomorphism $x\mapsto x\cap1$ with multiplication by $k$ on
$\BZ_n$. Hence $k\in H_3(\BZ_n)$ is realizable if and only if $k$
is a unit in $\BZ_n$. Moreover by 3.1~B, two $3$-manifolds $M_0$,
$M_1$ with $H_1\cong\BZ_n$ representing classes $k_0$, $k_1$ in
$H_3(\BZ_n)$ with respect to {\it some\/} identifications
$H_1(M_i)\cong\BZ_n$, will be surgery equivalent if and only if
there is an automorphism of $\BZ_n$ which induces an automorphism of
$H_3(\BZ_n)$ sending $k$ to $k_1$. Multiplication by (a unit) $m$
on $\BZ_n$ induces multiplication by $m^2$ on $H_3(\BZ_n)$ (see
Proposition~3 of \cite{Ru} and Theorem~29.5 of \cite{Co}).
Therefore we have derived the following.

\proclaim{Corollary 3.9} The set $\SS(\BZ_n)$ of surgery
equivalence classes of closed oriented $3$-manifolds with
$H_1\cong\BZ_n$ is in bijection with the set of equivalence classes
of units of $\BZ_n$, modulo squares of units, the correspondence
being given by the image of the fundamental class of the manifold
in $H_3(\BZ_n)\equiv\BZ_n$. The correspondence is also given by the
self-linking linking number of a generator of $H_1$ ($\la(1,1)=q$
and $q$ is viewed as an element of $\BZ_n$). The equivalence class
of $q\in\BZ_n$ contains the lens space $L(n,q)$.
\endproclaim

\pf{3.9} The first statement is proved above using 3.1~B. For the
second statement, use 3.1~D. Note that since $H_1$ is cyclic, the
cup products $H^1\op H^1\to H^2$ must vanish by anticommutativity
with any coefficients unless $n$ is even. Therefore if $n$ is odd,
part~a) of 3.1~D is vacuous. For even $n$ the triple cup product
form on a cyclic group is determined by the linking form [Tu2;
Theorem~1], so in any case we need only consider condition~c) of
3.1~D. Clearly the linking form $\la$ on a cyclic group is
determined by $\la(1,1)=\f an$ and $a\in\BZ_n$ is well-defined
modulo squares of units. Moreover two such forms $\la$ and $\la'$
are isomorphic if and only if $a\equiv a'$ modulo squares. For
$L(n,q)$,\ \ $\la(1,1)=\f qn$. \qed

\proclaim{Corollary 3.10} ($H_1\cong\BZ_p$, $p$ prime): Any
$3$-manifold with $H_1\cong\BZ_2$ is surgery equivalent to $\BR
P(3)$. For any odd prime there are precisely two surgery
equivalence class represented by $L(p,1)$ and $L(p,q)$ where
$q\not\equiv k^2\mod p$. Hence if $p\equiv3\mod 4$ then we may
take $q=-1$ and we see that $L(p,1)$ and $-L(p,1)$ are {\it not\/}
surgery equivalent. If $p\equiv1\mod 4$ then $-L(p,1)$ is surgery
equivalent to $L(p,1)$.
\endproclaim

More generally we see that:

\proclaim{Corollary 3.11} If $M_0$ and $M_1$ are
orientation-preserving homotopy equivalent then they are surgery
equivalent.
\endproclaim

\pf{3.11} This is immediate from 3.1~B.

\proclaim{Corollary 3.12} If $\pi_1(M_0)\cong\pi_1(M_1)$ is
abelian then $M_0$ is surgery equivalent to $M_1$ if and only if
$M_0$ is orientation-preserving homotopy equivalent to $M_1$.
\endproclaim

\pf{3.12} One implication follows from 3.11. Suppose $M_1$ is
surgery equivalent to $M_0$. By Theorem~4.2 there is cobordism $W$
from $M_0$ to $M_1$ built from $M_0\x[0,1]$ by adding two-handles
attached along curves being in $[\pi_1(M_0),\pi_1(M_0)]=0$. Hence
$W\simeq M_0\vee S^2\dots\vee S^2$, and there is a retraction
$r:W\to M_0$. The inclusion $M_1\to W$ followed by $r$ is a degree~1
map $M_1\to M_0$ inducing an isomorphism on $\pi_1$ and all
homology groups. We may assume the manifolds contain no fake
$3$-cells since these are irrelevant to the question of being
homotopy equivalent. Since $\pi_1$ is abelian it is not a
non-trivial free product so we may assume that
$\pi_2(M_0)=\pi_2(M_1)=0$ or that $M_0\cong M_1\cong S^1\x S^2$. In
the first case $\pi_1$ must be finite cyclic and then it is easy to
see that $f$ induces an isomorphism on $\pi_3$ by considering the
universal cover of $W$. Hence $f$ is a degree~1 homotopy
equivalence. \qed

\proclaim{Corollary 3.13} $L(n,q)$ is surgery equivalent to
$L(n,q')$ if and only if they are orientation-preserving homotopy
equivalent, that is if $qq'\equiv k^2\mod n$ for some unit $k$.
\endproclaim

\ex{3.14 ($H_1\cong\BZ\x\BZ_n$)} Since
$H_3(\BZ\x\BZ_n)\overset\pi\to\lra H_3(\BZ_n)$ is an isomorphism,
the surgery equivalence class depends only on the linking form.
From another point of view, since $H^1(\BZ\x\BZ_n;\BZ_m)$ is
generated by $2$-elements, the triple cup product forms vanish if
$m$ is odd and are determined by the linking form if $m$ is $2^r$.
Therefore each surgery equivalence class contains a representative
of the form $S^1\x S^2\#L(n,q)$ and the self-linking number of an
element of order $n$ in $H_1(M)$, $\la(1,1)=\f an$ viewed as an
element $a\in\BZ_n$ will distinguish the classes when viewed in the
group of units modulo squares (as in the case $H_1\cong\BZ_n$).

\ex{3.15} Some words of caution are in order. One must be careful in
applying 3.1. There exist $3$-manifolds which have isomorphic $H_1$,
linking forms and {\it integral\/} triple cup product forms but are
not surgery equivalent as detected by a $\BZ_p$ triple cup product
form. Namely, let $M_0$ be $\#^3_{i=1}L(5,1)$ and let $M_1$ be
$5/1$-surgery on each component of a Borromean Rings.

It is even possible that $M_0$ and $M_1$ have isomorphic linking
forms and isomorphic triple cup product forms with all coefficients,
yet {\it not\/} be surgery equivalent because the isomorphisms are
not induced by the same isomorphism $\phi$ on $H_1$! Consider the
manifolds in Figure~3.16. Since they have the same linking matrix
(expand the $5/2$ to a chain of integral surgeries if you like),
their linking forms are isomorphic to $\(1/5\)\op\(2/5\)$ on
$\BZ_5\x\BZ_5$. The triple cup product form on integral $H^1$ is
zero while the triple cup product forms on $H^1(\un{\ \ };\BZ_5)
\equiv(\BZ_5)^4$ are isomorphic by ``swapping the meridians
of the $5$ and $5/2$ knots.'' But these isomorphisms are
incompatible. We do not provide details.

\midinsert
\vspace{1.2in}
\botcaption{Figure 3.16}\endcaption
\endinsert

We can relate surgery equivalence to two other geometric equivalence
relations which have appeared in the literature. The first is
concerned with Dehn surgery on links; the second is concerned with
Heegard splittings and mapping class groups.

Recall that, given $M_0$, any closed, oriented $3$-manifold $M_1$
can be obtained by integral surgery on {\it some\/} framed link in
$M_0$. If the links are restricted what can be said? Recall that a
{\it boundary link\/} is a very special link with all linking
numbers zero, namely one whose components bound disjoint Seifert
surfaces. This makes sense in any $3$-manifold. The following is
mildly surprising.

\proclaim{Proposition 3.17} Let $M_0$, $M_1$ be oriented, connected
$3$-manifolds. Then the following are equivalent:
\roster
\item"A." $M_0$ is $2$-surgery equivalent to $M_1$
\item"B." There is a framed boundary link $L$ in $M_0$ such that
Dehn surgery (with framings $\pm1$) on $L$ yields $M_1$.
\item"C." $M_0$ and $M_1$ have isomorphic linking forms and triple
cup product forms (in the sense of {\rm3.1}).
\endroster
\endproclaim

Once again the result for homology spheres, where C is vacuous,
was known. The earliest reference we can find is S.V.~Matveev
\cite{Ma; Theorem A}. This result was reproved and used by
S.~Garoufalidis \cite{Ga}.

\pf{3.17} The equivalence of A and C is part of 3.1. B$\Rightarrow$A
is almost immediate. One need only note that the remaining
components of a boundary link remain null-homologous (use the same
Seifert surface) after performing Dehn surgery on some of its
components. Similarly the framings remain $\pm1$ because the
longitude remains the same. Thus we need only establish
A$\Rightarrow$B. By 2.8 we can assume that $M_1$ is the result of
$\pm1$ surgeries on an link in $M_0$ whose pairwise linking numbers
are zero. By induction, suppose $M_1$ is $\pm1$ surgery on a
null-homologous knot $K$ in $\SS(M_0,\{L_1,\dots,L_n\})$, the result
of $\pm1$-framed surgery on the boundary link $L$ in $M_0$, where
$\ell k(K,L_i)=0$ for all $i$. The following type of argument has
been used by others to prove the case of homology spheres. It serves
equally well in general. Let $L'\sbq\SS$ be the link consisting of
the cores of the surgery solid torii (so $\SS-L'\equiv M_0-L$). We
shall describe an isotopy of $K$ to $K'$ in $\SS$ (passing through
$L'$) such that $L\cup K'$ is a boundary link in $M_0$. Consider a
set $\SW=\{W_1,\dots,W_n\}$ of disjoint Seifert surfaces for $L$ in
$M_0-L$ whose boundaries are longitudes. Since we may actually
assume the homeomorphisms defining the Dehn surgeries carry
longitudes to longitudes, these surfaces can be extended by adding
annuli to $\wh\SW=\{\wh W_1,\dots,\wh W_n\}$, disjoint Seifert
surfaces in $\SS$ for the components of $L'$. Now $K$ bounds a
surface $V$ in $M_0-L\equiv\SS-L'$ because of the hypothesis on
linking numbers. Hence $\p\wh W_i\cap V=\phi$. $V$ has a
$1$-dimensional spine whose transverse intersections with $\wh\SW$
may be removed by isotopy merely by pushing over $\p\wh\SW$. This
isotopy extends to $V$. The resulting $K'=\p V'$ forms a boundary
link with $L$ since $V'\cap\SW\sbq V'\cap\wh\SW=\phi$. Note that, in
the presence of other components $K_2$, $K_3,\dots$ such that
$\ell k(K,K_i)=0$, the isotopy can be chosen so as to preserve that
latter fact. \qed

\bpage
The reader might find it interesting to compare this with Theorem~B
of \cite{Ma} which maintains that $M_0$ and $M_1$ have isomorphic
homology groups and linking forms if and only if $M_1$ can be
obtained from $M_0$ by surgery on a ``$T_0$-boundary link''
(recently re-introduced by Garoufalidis and Levine who used the term
``blink'' \cite{GL}). This was an improvement on a theorem of Hilden
who showed that any homology $3$-sphere can be obtained from $S^3$
by surgery on a blink \cite{H}.

Another way to describe $3$-manifolds is by their Heegard
splittings. Given $M_0$, by choosing a Heegard splitting
$M_0=H_1\cup_fH_2$, one can vary $f$ by composing with another
homeomorphism $g$ and in doing so change the three-manifold. One can
then ask if an arbitrary $M_1$ may be obtained in such a manner (it
{\it can\/} by Lickorish's theorem), or if it can be obtained using
only $g$ taken from some subgroup of the group of homeomorphisms. In
particular, if $\G$ is the mapping class group for the closed,
orientable surface of genus $g$ then we can consider the {\it
Torelli group\/} $\ST$, and the {\it Johnson group\/} $\SK$ of $\G$.
The Torelli group is the subgroup generated by homeomorphisms
inducing the identity map on $H_1$. The {\it Johnson subgroup\/}
$\SK\sbq\ST$ is the subgroup of $\G$ generated by Dehn twists along
simple closed curves in the surface which bound a subsurface
\cite{Jo1} \cite{Mo1} \cite{Mo2} \cite{GL}. Then, building on the
observations of \cite{GL} and our own work we have the following
generalization of the theorem of Morita for homology spheres
\cite{Mo1; Prop.~2.3}.

\proclaim{Theorem 3.18} Let $M_0$, $M_1$ be closed oriented
$3$-manifolds. Then the following are equivalent.
\roster
\item"A." $M_0$ and $M_1$ are $2$-surgery equivalent.
\item"B." There exist Heegard splittings $M_0=H_1\cup_fH_2$,
$M_1=H_1\cup_{\psi\circ f}H_2$ where $\psi\in\SK$.
\item"C." $M_0$ and $M_1$ have isomorphic linking forms and triple
cup product forms (as in {\rm3.1}).
\endroster
\endproclaim

\pf{3.18} The arguments in the first 3 paragraphs of section 2.3 of
\cite{GL}, although given for homology spheres, suffice to show that
3.18B is equivalent to 3.17B. \qed

\rem{3.19} We have shown that $\SK$ corresponds to boundary links
which, in turn, corresponds to homology surgery equivalence.
Strangely, the seemingly more natural Torelli group corresponds to
``blinks'' (see \cite{Ma} and \cite{GL}) which corresponds to
preserving only $H_1$ and the linking form. Combining the work of
\cite{Ma}, \cite{GL} and our present work yields an analagous theorem
to this effect, with ``homology surgery equivalence'' changed so
that the relation is generated by surgery on a $2$-component blink
rather than a knot.

\subhead{\bf\S4. Proofs of Theorem 3.1 and other Basic
Theorems}\endsubhead  In this section we prove Theorem 3.1. Several
major components of the proof are derived in much greater
generality so that they can be employed in later sections.

\bpage
\noindent$\BA\boldsymbol\Rightarrow\BB$: We will prove a more
general result which will be useful later.

\proclaim{Proposition 4.1} Suppose $M_1$ is obtained from $M_0$ by
$\pm1$ surgeries on a link $\{\g_1,\dots,\g_n\}$ where $\g_i\in
N\l\pi_1(M_0)\r$ and the meridians $\g'_i\in N\l\pi_1(M_1)\r$. Then
there exists an isomorphism
$\phi:\pi_1(M_1)/N\l\pi_1(M_1)\r\to\pi_1(M_0)/N\l\pi_1(M_0)\r$ such
that $(f_0)_*([M_0])=(f_1)_*([M_1])$ in
$H_3\(\pi_1(M_0)/N\l\pi_1(M_0)\r;\BZ\)$ where $f_1$ induces $\phi$
on $\pi_1$ and $f_0$ is the ``natural inclusion.'' Equivalently,
there exists an $f_1$, inducing an isomorphism $\phi$, such that
$(M_0,f_0)$ and $(M_1,f_1)$ are bordant over
$K\(\pi_1(M_0)/N\l\pi_1(M_0)\r,1\)$.
\endproclaim

\pf{4.1} Let $W$ be the usual cobordism from $M_0$ to $M_1$ obtained
by attaching $2$-handles to the $\g_i$ in $M_0\x\{1\}\subset
M_0\x[0,1]$. The curves $\g'_i$ are the attaching circles of the
dual $2$-handles attached to $M_1$. Hence the inclusions $j_0$,
$j_1$ induce isomorphisms $(j_0)_*$, $(j_1)_*$ on $\pi_1$ ``modulo
the $N$-subgroup.'' Let $\phi=(j_0)^{-1}_*\circ(j_1)_*$. Let
$\psi:\pi_1(W)/N\l\pi_1(W)\r\to\pi_1(M_0)/N\l\pi_1(M_0)\r$ be
$(j_0)^{-1}_*$. Then there are continuous maps $f_0$, $F$ and $f_1$
from $M_0$, $W$, $M_1$ respectively inducing the projection, $\psi$
and $\phi$ respectively on $\pi_1$ and such that $F$ extends $f_i$.
The result follows. \qed

\bpage
To see that 4.1 implies $A\Rightarrow B$, we apply 2.3 to reduce to
the case of $\pm1$ surgeries, we note that it suffices to prove
$A\Rightarrow B$ for the case of a single $\pm1$ surgery, then
appeal to 2.1 to see that the hypotheses of 4.1 are satisfies.

\bpage
\noindent$\BA\boldsymbol\Rightarrow\BD$: It suffices to consider
the case that $M_0$ and $M_1$ are cobordant via a single $2$-handle
attached with $\pm1$ framing along a circle $\g$ which is
null-homologous in $M_0$. Let $\phi_1=(j_0)^{-1}_*\circ(j_1)_*$ be
the isomorphism on $H_1$ induced by the inclusions. Then $\phi_1$
algebraically induces $\phi^1_n$ as in $C$ and the key point is to
observe that in this case $\phi^1_n$ equals
$j^*_1\circ(j^*_0)^{-1}$ where $j^*_i$ are the isomorphisms on
$H^1(\un{\ \ };\BZ_n)$ induced by the inclusions. Hence for any
$\a$, $\b$, $\g$ in $H^1(M_0;\BZ_n)$, we have
$$
\<\phi^1_n(\a)\cup\phi^1_n(\b)\cup\phi^1_n(\g), [M_1]\> =
\<(j^*_0)^{-1}(\a)\cup(j^*_0)^{-1}(\b)\cup(j^*_0)^{-1}(\g),
(j_1)_*([M_1])\>
$$
since $j_1$ is a continuous map. But $(j_0)_*([M_0])=(j_1)_*([M_1]$
for any coefficients so the above expression equals
$\<\a\cup\b\cup\g,[M_0]\>$ as desired. This shows condition a) of D
and C. Now we demonstrate that
$\la_1(\phi^{-1}_1x,\phi^{-1}_1y)=\la_0(x,y)$ for all torsion
classes $x$, $y$ in $H_1(M_0;\BZ)$. We may choose circles $\tl x$,
$\tl y$ in $M_0$ to represent these classes so that they are
disjoint from $\g$ and are in fact disjoint from a Seifert surface
$S$ for $\ell(\g)$. This is true because if $\tl x$ hits $S$, we
are free to isotope $\tl x$ ``through $\g$'' to achieve that the
algebraic number of such intersections is $0$. Then modify $S$ to
miss $\tl x$. Choose an integer $n$ and a $2$-chain $d$ of $M_0$
such that $\p d=n\tl x$ and such that $d$ meets $\tl y$ and $\g$
transversely. Then $\la_0(x,y)$ equals $\f1n$ times $\#(\tl y\cd
d)$. Now note that $\tl x$ and $\tl y$ are perfectly good
representatives of $\phi^{-1}_1(x)$ and $\phi^{-1}_1(y)$ since
$\phi_1$ is induced by the inclusions. We construct $d'$, a
$2$-chain in $M_1$ such that $\p d'=\p d=nx'$. For each
intersection of $d$ with $\g$, delete the $2$-disk $d\cap N(\g)$
and replace it by a copy of the annulus in the surgery solid torus
which expresses the fact that $\mu(\g)$ is isotopic to
$\pm\ell(\g)$ after surgery, and a copy of $\pm S$. This $2$-chain
$d'$ lies in $M_1$, has the same boundary as $d$ and $d'\cd\tl
y=d\cd\tl y$ since $\tl y\sbq M_0-N(\g)$ and $\tl y$ is disjoint
from $S$. Hence $\la_1\l\phi^{-1}_1(x),\phi^{-1}_1(y)\r=\la_0(x,y)$.
\qed

\bpage
\noindent$\BB\boldsymbol\Rightarrow\BA$: We will prove a
significantly broader result than is necessary in order to use it
in later sections.

\proclaim{Theorem 4.2} Suppose $M_0$ and $M_1$ are closed, oriented
$3$-manifolds. Suppose $N\trianglelefteq\pi_1(M_0)$ is the normal
closure of a finite number of elements and that $N$ is contained in
the commutator subgroup of $\pi_1(M_0)$ (or merely in
$\l\pi_1(M_0)\r^\BQ_2$ in the rational case). Suppose there exists
an epimorphism $\phi:\pi_1(M_1)\to\pi_1(M_0)/N$ which induces an
isomorphism on $H_1(\un{\ \ };\BZ)$ (or merely $H_1(\un{\ \ };\BQ)$
in the rational case) such that $(f_0)_*([M_0])=(f_1)_*([M_0])$ in
$H_3\l\pi_1(M_0)/N;\BZ\r$ (here $f_1$ is induced by $\phi$ as
usual, and $f_0$ induced by the inclusion into
$K\l\pi_1(M_0)/N,1)\r$. Then $(M_0,f_0)=(M_1,f_1)$ in $\Om_3\l
K(\pi_1(M_0)/N,1)\r$ via a $4$-manifold with only $2$-handles
$(\rel M_0)$ whose attaching circles lie in $N$ and whose linking
matrix with respect to $M_0$ is diagonal and invertible over $\BZ$
(merely invertible over $\BQ$ in the rational case). Consequently
$M_0$ is obtained from $M_1$ by $\pm1$ surgeries (integral surgeries
in the rational case) on a link $\{\g_1,\dots,\g_n\}$ such that
$[\g_i]\in N$. This link may be ordered in such a way that the
sequence of these surgeries exhibits that $M_1$ is $N$-surgery
related to $M_0$. Conversely $M_1$ is obtainable from $M_0$ in a
similar manner by surgery on a link $\{\g'_i\}$ and that
$\ker\phi$ is the normal subgroup of $\pi_1(M_1)$ generated by
$\{\g'_i\}$. Example~2.6 shows that this cannot, in general, be
strengthened.
\endproclaim

Before proving 4.2, we should check that it implies $B\Rightarrow
A$. Apply 4.2 with $N=\l\pi_1(M_0)\r_2$ to get the link
$\{\g_1,\dots,\g_n\}$. Apply 2.8. Hence $M_1$ is $2$-equivalent via
integral surgeries to $M_0$.

\pf{4.2} Let $X=K\l\pi_1(M_0)/N,1\r$ which we may think of as
constructed by adjoining cells to $M_0$. Then we have natural maps
$f_0:M_0\to X$ and $f_1:M_1\to X$ such that $(f_1)_*=\phi$. It is
well known that the map from $\Om_3(X)\to H_3(X)$ is an isomorphism
given by the image of the fundamental class. Thus the hypotheses
guarantee that there is a compact oriented $4$-manifold $W$ and a
map $F:W\to X$ such that $\p(W,F)=(M_1,f_1)\amalg(-M_0,f_0)$. Since
$F_*$ is necessarily surjective on $\pi_1$, we may perform surgery
on circles in $W$ and assume $F_*$ is an isomorphism.

Choose a handlebody structure of $W\rel M_0$ with no handles of
index~0 or~4. We may then proceed to ``trade'' $1$-handles for
$2$-handles as in [Ki2, pp\.~6--7, p\.~247]. This may also be
thought of as performing a surgery on the interior of $W$ along a
circle $c$ passing over the $1$-handle. Since $(f_0)_*$ is an
epimorphism on $\pi_1$, these circles may be altered by loops in
$M_0$ so that $F_*(c)=0$ (then in fact these loops are
null-homotopic) and hence the map $F$ extends to the ``new'' $W$.
Since $\phi_*=(f_1)_*$ is surjective we may trade all $3$-handles for
$2$-handles, by viewing then as $1$-handles attached to $M_1$.

Now let $V$ be the ``linking matrix'' of the attaching maps of the
$2$-handles $\rel M_0$. By this we mean the following. If
$\{\g_1,\dots,\g_m\}$ denote the attaching circles and
$\{\rho_1,\dots,\rho_m\}$ are the surgery circles on $\p N(\g_i)$
which are null-homotopic in $M_1$, then let $v_{ij}=\ell
k(\rho_i,\g_j)$. Since $\rho_i$, $\g_j$ are disjoint oriented
circles in $M_0$ which are null-homologous (torsion in the
$\BQ$-case) $v_{ij}$ is a well-defined integer (rational number in
the $\BQ$-case). We shall show that $V$ is invertible over $\BZ$
(respectively over $\BQ$).

We treat the $\BQ$ case first. Consider the long exact sequence in
rational homology for the pair $(M_0,M_0-\dot N(L))$ where
$L=\{\g_1,\dots,\g_m\}$. 
$$
\lra H_2\(M_0, M_0 - \dot N(L)\)\overset\p_*\to\lra H_1\(M_0 - \dot
N(L)\)\overset\pi_0\to\lra H_1(M_0)\lra0.
$$
Since the first term is $\BQ^m$ generated by the meridional disks 
we get $\BQ^m\overset i\to\lra H_1\l M_0-\dot N(L)\r
\overset\pi_0\to\lra H_1(M_0)\lra0$ where $i(e_i)=\mu_i=\mu(\g_i)$.
It is well-known that $i$ is injective when $\g_i$ are zero in
$H_1(\un{\ \ };\BQ)$. A splitting of $i$ is given by
$\phi(x)=\sum^m_{j=1}\ell k(x,\g_j)e_j$. Hence $\psi:H_1(M-L)\to
Q^m\op H_1(M_0)$ by $x\to\l\phi(x),\pi_0(x)\r$ is an isomorphism.
Note that $\psi(\rho_i)=\(\sum^m_{j=1}\ell k(\rho_i,\g_j)e_i,0\)$
since $\rho_i=0$ in $H_1(M_0;\BQ)$. Since $H_1(M_1)\cong
H_1(M_0-L)/\<\rho_i\>$, $H_1(M_1)\cong H_1(M_0)\op\BQ^{m-\rank V}$.
But $\phi$ is an isomorphism on $H_1(\un{\ \ };\BQ)$ so $\rank V=m$.

In the integral case we have the exact sequence in integral
homology:
$$
\BZ^m\overset i\to\lra H_1\(M_0 - \dot N(L)\)\overset\pi_0\to\lra
H_1(M_0)\lra0
$$
where $i$ is injective because it is injective with rational
coefficients. Here $i(e_i)=\mu_i$ and there is an isomorphism
$\phi:\ker\pi_0\to\BZ^m$ given by $\phi(x)=\sum^m_{j=1}\ell
k(x,\g_j)e_j$ such that $\phi\circ i=$~identity. In particular
$\phi(\rho_i)=\sum^m_{j=1}\ell k(\rho_i,\g_j)e_j
=\sum^m_{j=1}v_{ij}e_j$ so $i\(\sum^m_{j=1}v_{ij}e_j\)=\rho_i$.
Since $H_1(M_1;\BZ)\cong H_1\l M_0-\dot N(L)\r/\<\rho_i\>$, we see
that the cokernel of $\BZ^m\overset V\to\lra\BZ^m$ embeds in
$H_1(M_1;\BZ)$ via the map $i$. If this cokernel is non-zero then
there is a class $x\in\ker\pi_0$ which is {\it non-zero\/} under
$\pi_1:H_1\l M_0-\dot N(L)\r\lra H_1(M_1)$. This implies that
$\pi_1(x)$ is a non-trivial element in the kernel of the inclusion
map $H_1(M_1)\lra H_1(W)$ which is a contradiction. Hence $V$ is
invertible over $\BZ$.

Now that we have established that the linking matrix is invertible,
in the integral case, we appeal to the classification of symmetric
bilinear forms. We can change $W$ by adding a single $\pm1$ framed
$2$-handle attached along a trivial circle in order to assume $V$ is
an indefinite, odd form. Such a form has a $\BZ$-basis for which $V$
is diagonal with $\pm1$'s on the diagonal. This basis change can be
realized geometrically by handle slides (see [Ki2; Chapter~2]). Thus
$M_1$ is obtained from $M_0$ by $\pm1$ surgeries on a link
$\{\g_1,\dots,\g_n\}$ such that each
$[\g_i]\in N\l\pi_1(M_0)\r$ and $\ell k(\g_i,\g_j)=0$. Conversely
$M_0$ is obtainable from $\pm1$ surgery on the dual link
$\{\g'_1,\dots,\g'_n\}$ where, in general, all we know is that
$[\g'_i]\in\kernel\phi$ which in turn lies in
$[\pi_1(M_1),\pi_1(M_1)]$.

For the rational case we need the following Lemma, which may be well
known. The proof was suggested to me by Richard Stong.

\proclaim{Lemma 4.3} If $q:\BZ^n\x\BZ^n\lra\BQ$ is a symmetric,
bilinear, non-singular form then there is a basis
$e_1,\dots,e_n$ for $\BZ^n$ such that $q$ restricted to
$\<e_1,\dots,e_i\>$ is non-singular for each $i\le i\le n$.
\endproclaim

\prf First we claim that for any such form there is a basis such
that $q(e_1,e_1)\neq0$. For an arbitrary basis
$\{e_1,\dots,e_n\}$, if some $j$ has $q(e_j,e_j)\neq0$ then we are
done by re-ordering. If all $q(e_j,e_j)=0$ then, by
non-singularity, there is some $j$ such that $q(e_1,e_j)\neq0$.
Then the basis $\{e_1+e_j,e_2,\dots,e_n\}$ works.

Now we proceed by induction. Suppose we have a basis
$\{e_1,\dots,e_n\}$ such that $q|\<e_1,\dots,e_j\>$ is non-singular
for each $j<i$. We shall re-choose $\{e_i,\dots,e_n\}$ such that
$q|\<e_1,\dots,e_i\>$ is also non-singular. To do so write $q$ in
our basis as $q=\pmatrix  A  &B\\  B'  &C\endpmatrix$ where $A$ is
$(k-1)$ by $(k-1)$ and $B'$ is the transpose of $B$. We can make a
{\it rational\/} change of basis to replace this matrix by
$q^\land=\pmatrix  A  &O\\  O  &D\endpmatrix$ where $D=C-B'A^{-1}B$.
Since $q^\land$ is non-singular (since $q$ is), $D$ is a
non-singular matrix. As in the first step of our proof, there is an
integral invertible matrix $P$ such that the $(1,1)$ entry of
$P'DP$ is non-zero. Use this matrix $P$ to change our basis
$\{e_1,\dots,e_{i-1},e_i,\dots,e_n\}$ to
$\{e_1,\dots,e_{i-1},e'_i,\dots,e'_n\}$. In this new basis the
matrix for $q$ is obtained by conjugating the $q$ matrix by
$\pmatrix  I  &O\\  O  &P\endpmatrix$, which yields $q=\pmatrix  A 
&BP\\  P'B'  &P'CP\endpmatrix$. We claim that $q$ restricted to the
subspace spanned by $\{e_1,\dots,e_i\}$ is non-singular. To verify
this it suffices to apply the ``same'' {\it rational\/} change of
basis we used above but only to the first $i\x i$ submatrix. This
yields $\pmatrix  A  &O\\ O  &(P'DP)_{11}\endpmatrix$ where $11$
means the $(1,1)$ entry. Since this matrix is non-singular, the
original $q$ restricted to the span of $\{e_1,\dots,e_i\}$ is
non-singular. \qed

\bpage
Now apply 4.3 to the linking matrix $V$. The change of bases can be
achieved by re-ordering the handles and by ``handle slides'' in
$M_0$. Thus we may assume that $V_i$, the linking matrix of 
$\{\g_1,\dots,\g_i\}$ with respect to $M_0$, is non-singular for
each $1\le i\le n$. Let $M^*_i$ be the result of the surgeries on
$\{\g_1,\dots,\g_i\}$. We can assume by induction that the
cobordism $W^*$ from $M_0$ to $M^*_i$ is a product on $H_1$ modulo
torsion. Thus $[\g_{i+1}]$ is of finite order in $H_1(M^*_i;\BZ)$
since it is of finite order in $H_1(M_0;\BZ)$. By the argument of
the proof of 4.2, $H_1(M^*_{i+1};\BQ)\cong H_1(M_0;\BQ)$ since
$V_{i+1}$ is non-singular. But then $H_1(M^*_{i+1};\BQ)\cong
H_1(M^*_i;\BQ)$ so the surgery on $\g_{i+1}\sbq M^*_i$ is
non-longitudinal with respect to $M^*_i$. Thus $M_1$ is $N$-surgery
related to $M_0$.

This concludes the proof of Theorem~4.2. \qed

\bpage
We can now prove 2.9.

\pf{2.9} Suppose $M_0$ and $M_1$ are rationally $2$-surgery
equivalent. By 2.3 we may assume that there is a sequence
$M_0=X_0\to X_1\to\dots\to X_m=M_1$ where $X_{i+1}$ is obtained from
$X_i$ by a single integral non-longitudinal surgery on a rationally
null-homologous circle $\g_{i+1}$. We may assume $\{\g_i\}$ are
disjoint in $M_0$. Consider the induced cobordism $W$ from $M_0$ to
$X_i$ discussed in the proof of 2.3. The proof of 2.3 shows that $W$
is a product on $H_1$ modulo torsion. Since $[\g_{i+1}]$ is trivial
in $H_1(X_i;\BQ)$, it is also trivial in $H_1(M_0;\BQ)$. Moreover
the argument above in the proof of 4.2 shows that the linking matrix
of $\{\g_1,\dots,\g_m\}$ in $M_0$ is non-singular over $\BQ$. \qed

\bpage
We want to show that there exist maps $f_i:M_i\to X$ which induces
isomorphisms on the first integer homology group and such that
$(f_0)_*([M_0])=(f_1)_*([M_1])\in H_3(X)$. We may assume as before
that $f_0$ is the inclusion map. Here $X=K\l H_1(M_0),1\r$. Of
course, since $X$ is aspherical, there exists a map $f_1$ induced
by $\phi_1$.

First we note that it would suffice to show that
$\<f^*_0(k),[M_0]\>=\<f^*_1(k),[M_1]\>$ for {\it certain\/}
$h\in H^3(X;\BZ_n)$. For if $\a=(f_0)_*([M_0])-(f_1)_*([M_1])$ is
not zero in $H_3(X;\BZ)$ then there is an element of
$\hom(H_3(X);\BZ_n)$ which detects it, since $H_3(X)$ is a finitely
generated abelian group which has a element of order $n$ only if
$H_1(M_0)$ has an element of order $n$. More precisely, suppose
$H_1(M_0)\cong\BZ^m\x\BZ_{n_1}\x\dots\x\BZ_{n_k}$ where each $n_i$
is a prime power. Then the torsion-free summand of
$H_3\l H_1(M_0);\BZ\r$ is merely $H_3(\BZ^m;\BZ)\cong
H_3(S^1\x\dots\x S^1;\BZ)$. If $\a$ lies in this summand then it
can be detected by an element $h$ of the subgroup $H^3(\BZ^m;\BZ)$.
On the other hand if $\a$ is of finite order, then it can be
detected by some $h\in H^3(X;\BZ_n)$ where $n=p^r$ where $p^r$ is
the {\it maximal\/} order of all elements in $H_1(M_0)$ which have
order a power of $p$. This is true since $\BZ_{p^i}$ injects into
$\BZ_{p^r}$ if $r\ge i$. Therefore these are the only types of
elements $h$ we need consider.

Next we need to understand the cohomology rings $H^*(X;\BZ_n)$.

\proclaim{Proposition 4.4} Suppose $X$ is a finitely generated
abelian group and $n$ is the exponent of the $p$-torsion subgroup
of $X$ (elements of order $p^i$). Then the ring $H^3(X;\BZ_n)$ is
generated by elements of the form $\a\cup\b\cup\g$ and
$\a\cup\tau_*B(\g)$ where $\a$, $\b$, $\g\in H^1(X;\BZ_n)$, $B$ is
the Bockstein associated to $0\lra\BZ\overset n\to\lra
\BZ\overset\tau\to\lra\BZ_n\lra0$, and
$\tau_*:H^2(X;\BZ)\lra H^2(X;\BZ_n)$.
\endproclaim

Before proving 4.4, we finish the proof that $C\Rightarrow B$.
First consider the case that $h=\a\cup\b\cup\g$. Note that this
includes the case $n=0$ since $H^3(\BZ^m;\BZ)$ is generated by such
elements. Then we have (using C~b)),
$$
\split
\<f^*_0(h), [M_0]\>  &= \<f^*_0(\a)\cup f^*_0(\b)\cup f^*_0(\g),
[M_0]\>\\
&= \<\phi^1_n\circ f^*_0(\a)\cup\phi^1_n\circ f^*_0(\b)\cup
\phi^1_n\circ f^*_0(\g), [M_1]\>\\
&= \<f^*_1(h), [M_1]\>
\endsplit
$$
since $\phi^1_n\equiv f^*_1\circ(f^*_0)^{-1}$. Now suppose
$h=\a\cup\tau_* B(\g)$. $f^*_0\l\a\cup\tau_*B(\g)\r=
f^*_0(\a)\cup f^*_0\tau_* B(\g)=f^*_0(\a)\cup\tau_*Bf^*_0(\g)$
since $f_0$ is a continuous map. Hence, using condition~c of C we
have that:
$$
\split
\<f^*_0(h), [M_0]\>  &= \<\phi^1_n\circ f^*_0(\a)\cup\tau_*B
\(\phi^1_n\circ f^*_0(\g)\), [M_1]\>\\
&= \<f^*_1(\a)\cup\tau_*B\l f^*_1(\g)\r, [M_1]\>\\
&= \<f^*_1(h), [M_1]\>.
\endsplit
$$
Thus 4.4 will complete $C\Rightarrow B$.

\pf{4.4} First we need the following:

\proclaim{Lemma 4.5} Suppose $X$ and $Y$ are spaces whose homology
is finitely generated in each dimension. Suppose $n$ is a prime
power. Then the cohomology cross product induces an isomorphism:
$$
\T_n: \sum_{p+q=3}H^p(X;\BZ_n)\op_{\BZ_n}H^q(Y;\BZ_n)\lra
H^3(X\x Y;\BZ_n).
$$
\endproclaim

\pf{4.5} $\T_n$ is a monomorphism by [Mu; Theorem~61.6]. Since the
homology groups of $X$ and $Y$ are finitely generated the domain
and range of $\T_n$ are finite groups. Thus it will suffice to show
that they are {\it abstractly\/} isomorphic. First we list some
abbreviations: $x_p\equiv H_p(X;\BZ)$, $y_q\equiv H_q(Y;\BZ)$,
$x^t_p\equiv\BZ_n$-torsion subgroup of $x_p$,
$e^x_p=\ext(x_p;\BZ_n)$, $\ox\equiv\ox_\BZ$,
$\ox_n\equiv\ox_{\BZ_n}$. If $A$, $B$ are finitely generated
abelian groups, then the following are easily verified:
$\hom(A;\BZ_n)\cong A\ox\BZ_n$, $e^x_p=x^t_p\ox\BZ_n\cong x^t_p$,
$A*B\cong A^t*B^t$, $(A\ox\BZ_n)\ox_n(B\ox\BZ_n)\cong(A\ox
B)\ox\BZ_n$. Expanding $H^3(X\x Y;\BZ_n)$ using the Universal
Coefficient Theorem for cohomology and then the Kunneth Theorem for
homology and applying the above, we get $\bigoplus_{p+q=3}(x_p\op
y_q)\op(x^t_1*y^t_1)\op x^t_2\op y^t_2\op\ext(x_1\op y_1;\BZ_n)$
all tensored with $\BZ_n$. On the other hand the domain of $\T_n$
may be expanded as $\bigoplus_{p+q=3}\[\(x_p\op\BZ_n\op
x^t_{p-1}\)\op_n\(y_q\op\BZ_n\op y^t_{q-1}\)\]$. Expanding and
comparing terms shows that these expressions are isomorphic, using
the fact that $\ext(x_1\op y_1;\BZ_n)\op x^t_1\op y^t_1\cong
\(x^t_1\op_n(y_1\op\BZ_n)\)\op\((x_1\op\BZ_n)\op_n y^t_1\)$ which
is easily seen by expressing $x_1$, $y_1$ as direct sums of their
``torsion and torsion-free parts.''

Now we show that 4.5 implies 4.4. Since $X$ is a finitely-generated
abelian group, it is a product $\x^k_{i=1}X_i$ cyclic groups of
infinite or prime-power order. Now apply 4.5 inductively. Recall
that if $\a\in H^p(X;\BZ_n)$, $\b\in H^q(Y;\BZ_n)$ then
$\a\x\b=\pi^*_1(\a)\cup\pi^*_2(\b)$ where $\pi_i$ are the
projections to the factors. Then 4.5 implies that $H^3(X;\BZ_n)$ is
generated by elements of the form
$\pi^*(\a)\cup\pi^*(\b)\cup\pi^*(\g)$,
$\pi^*(\a)\cup\pi^*(\Delta)$, and $\pi^*(\G)$ where
$\a\in H^1(X_i;\BZ_n)$, $\b\in H^1(X_j;\BZ_n)$, $\g\in
H^1(X_k;\BZ_n)$, $\Delta\in H^2(X_s;\BZ_n)$, $\G\in H^3(X_m;\BZ_n)$.
The only cases where $H^2(X_s;\BZ_n)$ is non-zero are
$H^2(\BZ_{p^s};\BZ_{p^r})$ where $s\le r$ since $n$ is the
exponent. Consider the coefficient sequence
$0\lra\BZ_{p^r}\overset i\to\lra\BZ_{p^{2r}}\overset\pi\to\lra
\BZ_{p^r}\lra1$. Then the induced map
$H^1(\BZ_{p^s};\BZ_{p^{2r}}\overset\pi_*\to\lra
H^1(\BZ_{p^s};\BZ_{p^r})$ is zero since the composition
$\BZ_{p^s}\overset\phi\to\lra\BZ_{p^{2r}}\overset\pi\to\lra
\BZ_{p^r}$ is zero for any $\phi$. Hence the Bockstein
$\wt B:H^1(\BZ_{p^s};\BZ_{p^r})\lra H^2(\BZ_{p^s};\BZ_{p^r})$ is an
isomorphism and consequently $\Delta=\wt B(\g)$ for some
$\g\in H^1(X_s;\BZ_n)$ and $\pi^*(\a)\cup\pi^*(\Delta)$ is
$\pi^*(\a)\cup\wt B\pi^*(\g)$. But in this case
$\tau_*:H^2(\BZ_{p^s};\BZ)\lra H^2(\BZ_{p^s};\BZ_{p^r})$ is an
isomorphism and so one sees that $\wt B=\tau_*B$ as desired.

The only cases where $H^3(X_m;\BZ_n)$ is non-zero are
$H^3(\BZ_{p^m};\BZ_{p^r})$ where $m\le r$. We shall show
$\G=\a\cup\Delta$ reducing to the case above. Let $L=L(p^m,1)$ be
the $3$-dimensional lens space. The map $L\overset i\to\lra
K(\BZ_{p^m},1)$ can be constructed by adding cells of dimension~4
and higher to $L$ and thus induces isomorphisms on first and second
cohomology and a monomorphism on $H^3$. By Poincar\'e Duality for
$L$, $i^*(\G)=i^*(\a)\cup i^*(\Delta)$ for some $\a\in
H^1(\BZ_{p^m};\BZ_{p^r})$ and $\Delta\in H^2(\BZ_{p^m};\BZ_{p^r})$.
Hence $\G=\a\cup\Delta$ as claimed. This completes the verification
4.5$\Rightarrow$4.4. \qed

\bpage
\noindent$\BD\boldsymbol\Rightarrow\BC$: For brevity let $\phi_*$
denote the isomorphism $\phi_1$ and $\phi^*$ denote its induced
adjoint $\phi^1_n$. Let $i$ denote the inclusion $\BZ_n\lra\BQ/\BZ$
where $1\lra\f1n$. Let $T_i(M)$ denote the torsion subgroup of
$H_i(M;\BZ)$. Since a linking pairing $\la$ is non-singular,
$\la(-,a)$ is an isomorphism $T_1(M)\lra\hom(T_1(M);\BQ/\BZ)$. For
any $\g\in H^1(M;\BZ_n)$, $i\<\g,-\>\in\hom(T_1(M);\BQ/\BZ)$ and we
let $D(\g)$ be its inverse under the above isomorphism. Thus
$D:H^1(M;\BZ_n)\lra T_1(M)$ and $\la(D(\g),a)=i\<\g,a\>$ for all
$a\in T_1(M)$.

We claim that if $\g\in H^1(M_1;\BZ_n)$ then
$\phi_*D_0\phi^*(\g)=D_1(\g)$ where $D_0$, $D_1$ correspond to
$M_0$, $M_1$ respectively. For if $a\in T_1(M_1)$ then
$a=\phi_*(b)$ for some $b\in T_1(M_0)$. So
$\la_1(\phi_*D_0\phi^*(\g),a)=\la_1(\phi_*D_0\phi^*(\g),\phi_*(b))=
\la_0(D_0\phi^*(\g),b)$, by hypothesis c) 4.1 D. Continuing,
$\la_0(D_0\phi^*(\g),b)=i\<\phi^*(\g),b\>=i\<\g,\phi_*(b)\>=i\<\g,a\>$.
Hence $\phi_*D_0\phi^*(\g)=D_1(\g)$.

Next we claim that the Poincar\'e dual of $B(\g)$ is $D_1(\g)$
where $B:H^1(M_1;\BZ_n)\lra H^2(M_1;\BZ)$ is the Bockstein
associated to $0\lra\BZ\overset\cd n\to\lra\BZ\lra\BZ_n\lra0$. We
need to show $\la_1(B(\g)\cap[M_1],a)=i\<\g,a\>$ for each $a\in
T_1(M_1)$. Suppose $g$ is $1$-cochain with coefficients in $\BZ_n$
representing $\g$ and $\tl g$ is an integral $1$-cochain reducing
to $g$. Then $B(\g)$ is represented by $\f1n\d\tl g$, a $2$-cochain
with integral values. Note that $B(\g)$ is $n$-torsion so its
Poincar\'e dual lies in $T_1(M_1)$. If $\Sigma$ is a chain
representing the orientation class $[M_1]$ then the Poincar\'e dual
is represented by $\f1n\d\tl g\cap\Sigma=\f1n\p(\tl g\cap\Sigma)$.
Thus $\la_1(B(\g)\cap[M],a)$ is given by $\f1n\cd\#\l(\tl
g\cap\Sigma)\cd a'\r$ where $a'$ is a chain representing $a$ and
$\#$ is the number of signed intersection points modulo $n$. But
this is also a calculation of $i\<\g,a\>$, finishing the
verification of our second claim.

Now we can finish the proof that $D\Rightarrow C$. If $\a$,
$\g\in H^1(M_1;\BZ_n)$ then
$$
\split
i\<\a\cup\pi^*B\g, [M_1]\>  &= i\<\a,\pi^*B\g\cap [M_1]\>\\
&= i\<\a,\pi_*\l B(\g)\cap [M_1]\r\>\\
&= i\<\a,\pi_*D_1(\g)\>\\
&= i\<\a, D_1(\g)\>\\
&= \la_1\(D_1(\a), D_1(\g)\)
\endsplit
$$
and similarly
$$
\split
i\<\phi^*\a\cup\pi^*B\phi^*\g, [M_0]\>  &= \la_0\l D_0(\phi^*\a),
D_0(\phi^*(\g))\r\\
&= \la_1\l\phi_*D_0\phi^*\a, \phi_*D_0\phi^*\g\r
\endsplit
$$
by hypothesis D. By our first claim this equals $\la_1\l D_1(\a),
D_1(\g)\r$ as above. Since $i$ is injective, condition~b) of C is
established. \qed

\subhead{\bf\S5. Rational Homology Surgery Equivalence}\endsubhead 
In this chapter we address the question of when two $3$-manifolds
are related by a sequence of Dehn surgeries on rationally
null-homologous curves which preserve $H_1(\un{\ \ };\BQ)$. In the
language of \S2, this is case iv) where $N=G^\BQ_2=\{x\in
G\mid\exists n, x^n\in G_2\}$, or {\it rational $2$-surgery
equivalence\/}. Note that $G/G^\BQ_2=H_1(G)/T_1(G)$. By 2.2 this is
an equivalence relation and by 2.3 it is sufficient to consider
non-zero {\it integral\/} framings. We find that this is completely
controlled by the isomorphism class of the integral cup product form.
Beware that, because we restrict to rationally null-homologous
curves, a surgery which preserves $\b_1$ will necessarily preserve
$H_1/T_1$ and consequently $H^1(\un{\ \ };\BZ)$. Therefore it is not
possible to have, for example, $H_1(M_0)\cong\BZ$,
$H_1(M_1)\cong\BZ$ with the natural map between them being
``times~2.'' This would be possible to achieve by allowing certain
surgeries on curves in $M_0$ which are {\it essential\/} in
$H_1(M_0;\BZ)$ but not primitive. Hence rational $2$-equivalence is
NOT the relation generated by Dehn surgeries which preserve
$H_1(\un{\ \ };\BQ)$ but rather those which preserve $H_1(\un{\ \ };\BZ)/T_1$.

\proclaim{Theorem 5.1} Suppose $M_0$ and $M_1$ are closed, oriented
connected $3$-manifolds. The following 3 conditions are equivalent.
\roster
\item"$\BQ$A)" $M_0$ and $M_1$ are rationally $2$-surgery
equivalent; that is, each may be obtained from the other by a
sequence of non-longitudinal Dehn surgeries on circles which are
zero in $H_1(\un{\ \ };\BQ)$ (equivalently {\bf integral}
non-longitudinal surgeries.
\item"$\BQ$B)" There exists an isomorphism
$\phi_1:H_1(M_1)/T_1(M_1)\lra H_1(M_0)/T_1(M_0)$ such that
$(f_0)_*([M_0])=(f_1)_*([M_1])$ in $H_3\l
H_1(M_0)/T_1(M_0);\BZ\r\cong H_3\l(S^1)^{\b_1(M_0)};\BZ\r$ where
$f_0$ is induced by ``inclusion'' and $f_1$ is induced by $\phi_1$.
That is, $M_0$ and $M_1$ are bordant over $(S^1)^{\b_1(M_0)}$.
\item"$\BQ$C)" There exists an isomorphism $\phi_1$ as in the first
line of $\BQ B$ such that $\<\a\cup\b\cup\g, [M_0]\>=\mathbreak
\<\phi^1(\a)\cup\phi^1(\b)\cup\phi^1(\g),[M_1]\>$ for all $\a$,
$\b$, $\g\in H^1(M_0;\BZ)$ and $\phi^1$ is the adjoint (Hom-Dual)
of $\phi_1$. That is, the integral cup product forms of $M_0$ and
$M_1$ are isomorphic.
\endroster
\endproclaim

\pf{5.1} $\BQ\BA\boldsymbol\Rightarrow\BQ\BB$: By 2.3 with
$N=\l\pi_1(M_0)\r^\BQ_2$, we may reduce to the case of a single
integral non-longitudinal surgery on $\g\subset M_0$ such that
$[\g]\in N$. If $W$ is the cobordism corresponding to this surgery,
then $\b_1(W)=\b_1(M_0)$ and the inclusion map induces an
isomorphism modulo the $N$-subgroups (see last paragraph of proof
of 2.3), which, in the case at hand, means it induces an epimorphism
on $H_1$ modulo torsion. But by the symmetry of 2.1 and 2.2 the same
may be said of $M_1$. Hence we may let
$\phi_1=(j_0)^{-1}_*\circ(j_1)_*$ as in the proof of 4.1. Then, as in 4.1, there are continuous maps
$f_0$, $F$ and $f_1$ from $M_0$, $W$, $M_1$ respectively, to
$K(\BZ^{\b_1(M_0)},1)$ inducing the obvious maps on $\pi_1$ and the
result follows. \qed

\newpage
\noindent$\BQ\BA\boldsymbol\Rightarrow\BQ\BC$: Note that the
inclusion maps $j_0$, $j_1$ as defined above induce isomorphisms on
$H^1(\un{\ \ };\BZ)$ since $\hom(H_1;\BZ)\cong\hom(H_1/T_1;\BZ)$.
The proof is now the same as the first part of the proof of
$A\Rightarrow D$ in \S4.

\bpage
\noindent$\BQ\BC\boldsymbol\Rightarrow\BQ\BB$: Note that all the
maps $j^*_0$, $j^*_1$, $\phi^1$, $f^*_0$, $f^*_1$ are isomorphisms
on $H^1(\un{\ \ };\BZ)$. Since the cohomology ring of $K(\BZ^m,1)$
is well-known to be generated by triple cup products, the easy part
of the proof of $C\Rightarrow B$ in \S4 applies word for word. \qed

\bpage
\noindent$\BQ\BB\boldsymbol\Rightarrow\BQ\BA$: Apply 4.2 with
$N=\l\pi_1(M_0)\r^\BQ_2$. \qed

\bpage
Let us denote by $\SS^\BQ_m$ the set of rational homology surgery
equivalence classes of closed oriented $3$-manifolds with $\b_1=m$.
By 5.1, if $m<3$ then $\SS^\BQ_m$ contains a single element.

\proclaim{Corollary 5.2} If $m<3$ any 2 closed, oriented
$3$-manifolds with identical first Betti number $m$ are rational
homology surgery equivalent.
\endproclaim

\proclaim{Corollary 5.3} There is a bijection
$\SS^\BQ_m\to\La^3(\BZ^m)/\gl_m(\BZ)$ given by the integral triple
cup product form. Hence if $M_0$, $M_1$ have torsion-free homology
groups then they are rational homology surgery equivalent if and
only if they are integral homology surgery equivalent.
\endproclaim

\pf{5.3} See the proof of 3.5 and use Sullivan's work [Su].

\ex{5.4} It is possible for $M_0$ and $M_1$ to be rational homology
surgery equivalent, have isomorphic first homology and linking forms,
yet not be integral homology surgery equivalent (see Example~3.15).

\proclaim{Corollary 5.5} For any closed, oriented $3$-manifold $M$,
$M$ is rational homology surgery equivalent to $-M$.
\endproclaim

\subhead{\bf\S6. Surgery Equivalence Preserving Lower Central Series
Quotients}\endsubhead
We have seen that the relation generated by $\pm1/n$ Dehn surgery on
circles which lie in $(\pi_1(M))_k$ is an equivalence relation which
we called $k$-surgery equivalence. The equivalence relation
generated by non-longitudinal surgeries on circles lying in the
$k^\supth$ term of the rational lower central series, we call
rational $k$-surgery equivalence. Just as $2$-equivalence was
controlled by $G/G_2$ and the cup products (and linking form), we
shall see that $k$-equivalence is controlled by $G/G_k$ and higher
Massey products (and the linking form). We only attempt a complete
algebraic characterization of $k$-surgery equivalence to ``the zero
element,'' i.e. $\#^m_{i=1}S^1\x S^2$. Here we see that
$k$-equivalence is controlled by Massey products of length less than
$2k-1$, or equivalently by the isomorphism class of $G/G_{2k-1}$. A
similar characterization for the general case is made difficult by
the ill-definedness of Massey products and our ignorance of $H_3$ of
torsion-free nilpotent groups. It may well be, however, that there
is sufficient information in the literature to complete the general
characterization.

\proclaim{Theorem 6.1} Suppose $M_0$ and $M_1$ are closed, oriented,
connected $3$-manifolds. The following are equivalent.
\roster
\item"A)" $M_0$ and $M_1$ are $k$-surgery equivalent.
\item"B)" There exists an isomorphism
$\phi:\pi_1(M_1)/\l\pi_1(M_1)\r_k\lra\pi_1(M_0)/\l\pi_1(M_0)\r_k$
such that $\phi_*([M_1]^k)=[M_0]^k$ where $[M_i]^k$ means the image,
in $H_3\l\pi_1(M_i)/\l\pi_1(M_i)\r_k;\BZ\r$, of the fundamental class
of $M_i$ under some map $f_i:M_i\to
K\l\pi_1(M_i)/\l\pi_1(M_i)\r_k,1\r$ inducing the obvious quotient on
$\pi_1$.
\endroster
\endproclaim

\proclaim{Corollary 6.2} The set of $k$-surgery equivalence classes
of closed, oriented $3$-manifolds $M$ with
$\pi_1(M)/\l\pi_1(M)\r_k\cong G$ is in bijection with the subset of
$H_3(G)/\aut(G)$ consisting of those elements which are
``realizable'', that is which can arise as $[M]^k$ for {\bf some}
closed $3$-manifold. The correspondence is given by the fundmental
class (see {\rm[Tu1]} for an analysis of this realizable set.
\endproclaim

\pf{6.1} $A\Rightarrow B$ is implied by Propositions 2.3, 2.1 and
4.1.\newline
$B\Rightarrow A$ is implied by Theorem 4.2. \qed

\pf{6.2} Merely note that $[M_i]^k$ is only well-defined up to the
action of $\aut(G_i)$ on $H_3(G_i)$.

\proclaim{Theorem 6.3} Suppose $M_0$ and $M_1$ are closed, oriented,
connected $3$-manifolds. The following are equivalent.
\roster
\item"$\BQ$A)" $M_0$ and $M_1$ are rationally $k$-surgery equivalent.
\item"$\BQ$B)" Same condition as 6.1B with rational lower central
series replacing the integral one.
\endroster
\endproclaim

\proclaim{Corollary 6.4} The set of rational $k$-surgery equivalence
classes of closed, oriented $3$-manifolds with
$\pi_1(M)/\l\pi_1(M)\r^\BQ_k\cong G$ is in bijection with the subset
of $H_3(G)/\aut(G)$ corresponding to realizable classes (see
\cite{Tu1}).
\endproclaim

\pf{6.4} The argument for $\BQ A\Rightarrow\BQ B$ in the proof of
5.1 works here (using that finitely generated nilpotent groups are
Hopfian). For $\BQ B\Rightarrow\BQ A$ apply 4.2 with
$N=\l\pi_1(M_0)\r^\BQ_k$. \qed

\proclaim{Corollary 6.5} Any two $3$-manifolds with $\b_1=0$ (or any
two with $\b_1=1$) are rationally $k$-surgery equivalent for each
$k$.
\endproclaim

\pf{6.5} Suppose $\b_1(M_0)=1$. Then the epimorphism
$\pi_1(M_0)\twoheadrightarrow\BZ$ induces isomorphisms modulo any
term of the rational lower central series \cite{St}, so
$\pi_1(M_0)/\l\pi_1(M_0)\r^\BQ_k$ is $\BZ$. Since $H_3(\BZ)=0$, the
result follows from 6.4. \qed

\midinsert
\vspace{1in}
\botcaption{Figure 6.7}\endcaption
\endinsert

\ex{6.6} Let $M_0=S^1\x S^2\#S^1\x S^2$ and let $M_1$ be the
manifold obtained by $0$-surgery on each component of a Whitehead
link in $S^3$, as shown by the solid lines in Figure~6.7a.
Performing a $+1$ surgery on the dashed circle $\g$ in Figure~6.7b
transforms $M_0$ to $M_1$, and $\g$ is clearly null-homologous in
$M_0$ so $M_0$ and $M_1$ are $2$-surgery equivalent (as they must be
since $H_3(\BZ\x\BZ)=0$). However $\g\notin\l\pi_1(M_0)\r_3$ and so
it is not clear whether or not $M_0$ and $M_1$ are $3$-surgery
equivalent. In this case $M_1$ is known to be the ``Heisenberg
manifold'' (Euler class $\pm1$ circle bundle over the torus) whose
fundamental group is $F/F_3$ where $F$ is the free group on
$\{x,y\}$. Hence
$\pi_1(M_0)/\l\pi_1(M_0)\r_3\cong\pi_1(M_1)/\l\pi_1(M_1)\r_3\cong
F/F_3$. But 6.1B is not satisfied since $[M_0]$ represents the
trivial class in $H_3(F/F_3)$ since it factors through $H_3(F)$,
whereas $M_1$ is a $K(F/F_3,1)$ and so $[M_1]$ represents a
generator of $H_3(F/F_3)\cong\BZ$. Thus the manifolds are not
$3$-surgery equivalent (nor rationally $3$-surgery equivalent).

We shall now show that $k$-surgery equivalence is related to higher
order Massey products and that this is the correct generalization
of the triple cup product form. However, Massey products may not be
uniquely defined and this makes statements of results difficult. For
this reason we shall restrict our focus to situations where the
Massey products are uniquely defined. In general if $M_0$ and $M_1$
are $k$-surgery equivalent then their lower central series quotients
$G/G_j$ are isomorphic for $1\le j\le k$ and this is known to entail
a ``correspondence'' between order $k-1$ Massey products with any
abelian coefficients \cite{Dw; Corollary~2.7}. We shall state this
only for a restricted case.

\proclaim{Proposition 6.8} Suppose $M_0$ and $M_1$ are $k$-surgery
equivalent, and that all $j^\supth$ order Massey products vanish for
$M_0$ for $2\le j\le(k-2)$. Then there is an isomorphism
$\phi:H_1(M_1)\to H_1(M_0)$ such that for all abelian groups $A$ and
$\a_i\in H^1(M_0;A)$, $1\le i\le k$,
$\<\a_1\cup\<\a_2,\dots,\a_k\>,[M_0]\>=\<\phi^*\a_1\cup
\<\phi^*\a_2,\dots,\phi^*\a_k\>,[M_1]\>$ where $\<\a_2,\dots,\a_k\>$
is the Massey product in $H^2(M_0;A)$. In fact,
$\pi_1(M_0)/\l\pi_1(M_0)\r_k\cong F/F_k$ ($F$ a free group) if and
only if $H_1(M_0)$ is torsion free and all Massey products of length
less than $k$ vanish for $M_0$. More precisely, the latter
conditions imply that any given isomorphism
$\pi_1(M_0)/\l\pi_1(M_0)\r_{k-1}\cong F/F_{k-1}$ can be extended.
\endproclaim

\midinsert
\vspace{1.2in}
\botcaption{Figure 6.10}\endcaption
\endinsert

\ex{6.9} Let $M_0=\#^4_{i=1}S^1\x S^2$ and let $M_1$ be the
manifold shown in Figure~6.10a\. as $0$-surgery on the $4$-component
link $L$. Then $M_0$ is $2$-surgery equivalent to $M_1$ as shown by
the 3 circles labelled $\pm1$ in 6.10b.. Note these circles lie in
the $G_2-G_3$ where $G=\pi_1(M_0)$. It is known that the link in
6.10a\. has $\ov\mu(1234)=\pm1$, so for $M_0$ all Massey products
vanish, but for $M_1$, $\<x_1\cup\<x_2,x_3,x_4\>, [M_1]\>=\pm1$
where $x_i$ are the Hom-duals of the meridians \cite{Ko; Theorem~3}.
Hence $M_0$ and $M_1$ are not $4$ surgery equivalent, by 6.8. In
fact they are not even $3$-surgery equivalent but 6.8 is too weak to
show this.

\pf{6.8} Apply \cite{Dw; Corollary 2.7} to the maps $j_0:M_0\to W$
and $j_1:M_1\to W$ where $W$ is the cobordism over
$\pi_1(M_0)/\l\pi_1(M_0)\r_k$ guaranteed by 6.1. Then use naturality
of Massey products and $(j_0)_*([M_0])=(j_1)_*([M_1])$ to get the
first claimed result. Details are left for the reader.

We consider the last claim of 6.8. Suppose
$\pi_1(M_0)/\l\pi_1(M_0)\r_k\cong F/F_k$, $k\ge2$. Then there is a
map $f:M_0\to K(F/F_k,1)$ which induces an isomorphism on first
cohomology with any abelian coefficients. It is known that $F/F_k$
has vanishing Massey products of length less than $k$ and that
$H^2(F/F_k)$ is generated by Massey products of length $k$ \cite{O2;
Lemma~16}. By naturality,
$f^*\l\<x_1,\dots,x_m\>\r\sbq\<f^*x_1,\dots,f^*x_m\>$. Hence, if
$m<k$ then $0\in\<f^*x_1,\dots,f^*x_m\>$. But the first non-vanishing
level of Massey products are uniquely defined, so by induction,
$0=\<f^*x_1,\dots,f^*x_m\>$ for $m<k$. This implies that all Massey
products of length less than $k$ vanish for $M_0$.

Now suppose $H_1(M_0)$ is torsion-free and all Massey products of
length less than $k$ vanish for $M_0$. By induction we can assume
$\pi_1(M_0)/\l\pi_1(M_0)\r_{k-1}\cong F/F_{k-1}$ so there is a map
$g:F\to G$ ($\pi_1(M_0)=G$) inducing this isomorphism. It would then
suffice to prove
$$
\def\pp{@>\pretend{}\haswidth{\text{eeeee}}>>}
\def\p{@>\pretend \pi_*\haswidth{\text{eeeee}}>>}
\split
&0\pp H_2(F/F_{k-1})\pp F_{k-1}/F_k\pp0\\
&\Bigl\downarrow\hskip50pt\cong\Bigl\downarrow g_*
\hskip50pt\Bigl\downarrow g_*\\
H_2&(G)\p H_2(G/G_{k-1})\pp G_{k-1}/G_k\pp0
\endsplit
$$
that $F_{k-1}/F_{k-2}\cong G_{k-1}/G_k$ and the diagram above shows
that this is equivalent to showing $\pi_*$ is the zero map. Since
$H_2(M_0)$ maps surjectively to $H_2(G)$ it suffices to show
$H_2(M_0)\overset\pi_*\to\lra H_2(G/G_{k-1})$ is the zero map. But
$H^2(G/G_{k-1})\cong H^2(F/F_{k-1})$ is generated by
$\<x_1,\dots,x_{k-1}\>$ so
$\pi^*\<x_1,\dots,x_{k-1}\>=\<\pi^*x_1,\dots,\pi^*x_{k-1}\>=0$. Thus
$\pi^*$ and $\pi_*$ are zero. Note that we have actually proved that
any isomorphism $g_*:G/G_{k-1}\lra F/F_{k-1}$ can be extended to
$g_*:G/G_k\lra F/F_k$. \qed

\bpage
Example 6.9 indicates that 6.8 is too weak. Indeed, 6.1B should be
seen as two conditions, and 6.8 says that the first of these
conditions controls Massey products of lengths up to $k-1$. We shall
see that the second conditions controls lengths up to $2k-2$. This
is analagous to the Cochran-Orr conjecture that a link in $S^3$ is
``null $k$-cobordant'' if and only if its Milnor $\ov\mu$-invariants
of length $j$, $1\le j\le2k$ are zero. This has been positively
resolved by X.S.~Lin and Orr-Igusa \cite{L} \cite{IO}. Based on
techniques of the latter, we shall now discuss an algebraic
characterization of $k$-surgery equivalence to $\#S^1\x S^2$.

\proclaim{Theorem 6.10} For any integer $k\ge2$, $M$ is $k$-surgery
equivalent to $\#^m_{i=1}S^1\x S^2$ if and only if
$H_1(M)\cong\BZ^m$ and all Massey products of order {\bf less}
than $2k-1$ vanish for $M$. In particular if $M$ is zero surgery on
an $m$-component link $L$ in a homology $3$-sphere then $M$ is
$k$-surgery equivalent to $\#^m_{i=1}S^1\x S^2$ if and only if
Milnor's $\ov\mu$-invariants of length {\bf less} than $2k$
vanish for $L$.
\endproclaim

\pf{6.10} Let $G=\pi_1(M)$. Suppose $H_1(M)\cong\BZ^m$ and all
Massey products of length less than $2k-1$ vanish for $M$. By the
proof of the last part of 6.8, any isomorphism $g_*:G/G_k\lra F/F_k$
extends to an isomorphism $h_*:G/G_{2k-1}\lra F/F_{2k-1}$. It
follows that ``the'' natural map $M\overset\pi\to\lra
K(G/G_k,1)\overset g\to\lra K(F/F_k,1)$ factors through
$K(F/F_{2k-1},1)$ and thus that $[M]^k=0$ in $H_3(F/F_k;\BZ)$ since
Igusa and Orr have shown that the map $H_3(F/F_{2k-1})\lra
H_3(F/F_k)$ is zero \cite{IO}. Hence, by 6.1, $M$ is $k$-surgery
equivalent to $\#^m_{i=1}S^1\x S^2$.

Now suppose $M$ is $k$-surgery equivalent to $\#^m_{i=1}S^1\x S^2$.
Then $[M]^k=0$ in $H_3(F/F_k)$. First we show that this implies that
$G/G_{k+1}\cong F/F_{k+1}$. This follows from this more general
result.

\proclaim{Lemma 6.11} Suppose $\pi_1(M_0)=G_0$, $\pi_1(M_1)=G_1$,
$G_0/(G_0)_k\cong F/F_k$, $G_1/(G_1)_{k+1}\cong F/F_{k+1}$. If $M_0$
is $k$-surgery equivalent to $M_1$ then $G_0/(G_0)_{k+1}\cong
F/F_{k+1}$ by an isomorphism extending $f$.
\endproclaim

\pf{6.11} Consider a cobordism $W$ from $M_0$ to $M_1$ which
contains only $2$-handles and is an $F/F_k$-cobordism (see 4.2). Let
$F$, $F_0$, $F_1$ denote the maps to $K(F/F_k,1)$ from $W$, $M_0$,
$M_1$ respectively. For any cohomology classes
$\{x_1,\dots,x_k\}\subset H^1(M)$, choose $y_i\in H^1(F/F_k)$ so
$(F_0)^*(y_i)=x_i$. Then
$\<x_1,\dots,x_k\>=\<F^*_0y_1,\dots,F^*_0y_k\>=
j^*_0\<F^*y_1,\dots,F^*y_k\>$, where the Massey products are
uniquely defined since products of lesser length vanish since all
spaces have $G/G_k\cong F/F_k$ (see 6.9). By 6.8, it will suffice to
show $\<F^*y_1,\dots,F^*y_k\>=0$. Certainly
$j^*_1\<F^*y_1,\dots,F^*y_k\>=0$ since all Massey products of length
$k$ vanish for $M_1$. To finish we will show that $j^*_1:H^2(W)\lra
H^2(M_1)$ is injective on the image of $H^2(F/F_k)$. Recall that we
can assume that $W$ is built from $M_1\x[0,1]$ by adding $2$-handles
whose attaching circles lie in $F_k$. Thus $H_2(W)$ splits as
$H_2(M_1)\op H_2(W,M_1)$ where the latter is a free abelian group
generated by the cores of the $2$-handles ``capped-off'' by surfaces
in $M_1$. But these latter classes clearly become spherical in
$K(F/F_k,1)$ and hence vanish in $H_2(F/F_k,1)$. Considering the
dual splitting of $H^2(W)$, this implies that the image of
$H^2(F/F_k)$ lies in the summand $H^2(M_1)\hookrightarrow H^2(W)$.
It follows that $j^*_1$ is injective on this image. \qed

\bpage
Suppose $\pi_1(M)\cong G$ and $G/G_k\cong F/F_k$ for some free group
$F$. Suppose also that $f$ is a specific such isomorphism. Then we
can define $\th_k(M,f)\in H_3(F/F_k)$ to be the image of the
fundamental class under the map $M\lra K(F/F_k,1)$ induced by $f$.

\proclaim{Lemma 6.12}
$\th_k(M,f)\in\img\l H_3(F/F_{k+1})\overset\pi_*\to\lra H_3(F/F_k)\r$
if and only if there is some isomorphism $\tl f:G/G_{k+1}\lra
F/F_{k+1}$ extending $f$ such that $\pi_*\l\th_{k+1}(M,\tl
f)\r=\th_k(M,f)$.
\endproclaim

\pf{6.12} Suppose $\th_k(M,f)=\pi_*(x)$ $x\in H_3(F/F_{k+1})$. By
\cite{O1; Theorem~4}, there exists a $3$-manifold $M_1$ with
$\pi_1(M_1)\cong G_1$ and $G_1/(G_1)_k\overset g\to\cong F/F_{k+1}$
such that $\th_{k+1}(M_1,g)=x$. Since
$\pi_*(x)=\th_k(M_1,\pi\circ g)=\th_k(M,f)$, it follows from 6.1
that $(M,f)$ and $(M_1,\pi\circ g)$ are cobordant over $F/F_k$ and
are $k$-surgery equivalent. From 6.11 we can conclude that there is
an isomorphism $\tl f$ as desired. Since $\tl f$ extends $f$,
$\pi_*\l\th_{k+1}(M,\tl f)\r=\th_k(M,f)$, by definition. The other
implication of 6.12 is immediate. \qed

\bpage
To finish the proof of 6.10 we need the following theorem from
\cite{IO}.

\proclaim{Theorem 6.13 (Igusa-Orr [IO])} Suppose $F$ is the
free group of rank $m$. Let $N_i=\rank H_2(F/F_i)$. Then
$H_3(F/F_k;\BZ)$ is $\bigoplus^{2k-2}_{i=k}\BZ^{mN_i-N_{i+1}}$. If we
define the {\bf weight} of the summand $\BZ^{mN_i-N_{i+1}}$ to be
$i+1$ then the natural maps $H_3(F/F_{k+1})\lra H_3(F/F_k)$ have the
property that they preserve weight whenever possible and map each
weight summand injectively and by the zero map otherwise.
\endproclaim

\proclaim{Corollary 6.14} The map $H_3(F/F_{2k-1})\lra H_3(F/F_k)$
is zero. Any element in the kernel of $H_3(F/F_{k+j})\lra
H_3(F/F_k)$, $j\le k-1$, lies in the image of $H_3(F/F_{2k-1})\lra
H_3(F/F_{k+j})$.
\endproclaim

\pf{6.14} An element of the kernel of the above map has weight
greater than $2k-1$ and at most $2k+2j-1$ i.e\. at most $4k-3$,
which is precisely the range of weights possible for an element of
$H_3(F/F_{2k-1})$. \qed

\bpage
Finally, if $\th_k(M,f)=0$ for some $f$ then assume by induction
(use 6.12 for the start of induction) that $f$ extends to $\tl
f:G/G_{k+j}\lra F/F_k$ so that $\th_{k+j}(M,\tl f)$ is defined and
is a lift of $\th_k(M,f)$. But then as long as $j\le k-1$, 6.14
guarantees that $\th_{k+j}(M,\tl f)\in\img\l H_3(F/F_{k+j+1})\lra
H_3(F/F_{k+j})\r$ and 6.12 then implies $\th_{k+j+1}$ is defined.
Hence we can lift to $\th_{2k-1}$ and $G/G_{2k-1}\cong F/F_{2k-1}$
implying that all Massey products of length less than $2k-1$ are
zero for $M$. This concludes the proof of the first sentence of
6.10. The second statement follows immediately from \cite{Ko;
Theorem~2} relating Massey products of $0$-surgery on a link to
$\ov\mu$-invariants of a link. \qed

\bpage
It follows from 6.10 that the two manifolds in Example 6.9 are not
$3$-surgery equivalent.

\ex{6.15} The example $M_1$ in Figure 6.7a (the Heisenberg manifold)
is $2$-surgery equivalent to $S^1\x S^2\#S^1\x S^2$ since the first
non-zero Milnor invariant of the Whitehead link, $\ov\mu(1122)$, is
of length~$4$, but for the same reason the Whitehead link is {\bf
NOT} null $2$-cobordant. In general, if a link is null $k$-cobordant
then its $0$-surgery is $k$-surgery equivalent to $\#S^1\x S^2$ but
the converse is false, requiring in addition that the
$\ov\mu$-invariants of length $2k$ vanish.

\vskip.7cm
\noindent Tim Cochran\newline
Mathematics Department

\noindent Rice University

\noindent 6100 Main Street

\noindent Houston, TX 77005--1892

\noindent cochran\@rice.edu

\noindent Related papers available also at
http://math.rice.edu/$\sim$cochran/

\end

For $2$-equivalence, doing surgery on such a
link is {\it equivalent\/} to doing surgeries on successive
$\p,1$-framed, null-homologous knots. To see this let $Y_i$ be the
result of surgery on $\{\g_1,\dots,\g_i\}$. Since $\ell
k(\g_i,\g_j)=0$, there is a Seifert surface for $\ell(\g_i)$ in
$M_0-\bigcup^n_{j=1}\g_j$. Thus $\g_i$ is null-homologous in $Y_i$
and the surgery framing is still $\pm1$ with respect to $Y_i$.

\bpage
Now, in the rational case we argue that any such surgeries are
non-longitudinal. From the above discussion we see that
$\b_1(M_0-L)=m+\b_1(M_0)$ and so if $[\rho_i]=0$ in $H_1(M_0-L;Q)$
then $\b_1(M_1)>\b_1(M_0)$ since $H_1(M_1)\cong
H_1(M_0-L)/\<\rho_i\>$.